\documentclass{cmslatex}
\usepackage[paperwidth=7in, paperheight=10in, margin=.875in]{geometry}
 \usepackage[backref,colorlinks,linkcolor=red,anchorcolor=green,citecolor=blue]{hyperref}
 \usepackage[hyphenbreaks]{breakurl}
\usepackage{amsfonts,amssymb}
\usepackage{amsmath}
\usepackage{graphicx}
\usepackage{cite}
\usepackage{enumerate}
\renewcommand{\theequation}{\arabic{section}.\arabic{equation}}

\usepackage{algorithm}
\usepackage{algorithmic}
\DeclareMathAlphabet{\mathpzc}{OT1}{pzc}{m}{it}

\allowdisplaybreaks

\usepackage{amsfonts}
\usepackage{graphicx}

\usepackage{indentfirst}
\usepackage{mathrsfs}
\usepackage{graphicx}
\usepackage{hyperref}
\usepackage{thmtools}
\usepackage{csquotes}
\usepackage{verbatim}
\usepackage{bookmark}
\usepackage{booktabs}
\usepackage{url}
\usepackage{hyperref}
\usepackage{enumerate}
\usepackage{dsfont}
\usepackage{xcolor}
\usepackage{mathtools}

\DeclareMathAlphabet{\mathpzc}{OT1}{pzc}{m}{it}

\newcommand{\thmref}[1]{Theorem~\ref{#1}}

\newcommand{\propref}[1]{Proposition~\ref{#1}}
\newcommand{\lemref}[1]{Lemma~\ref{#1}}
\newcommand{\secref}[1]{Sec.~\ref{#1}}

\newcommand{\ie}{{i.e.}}
\newcommand{\eg}{{e.g.}}

\makeatletter
\newcommand*{\rom}[1]{\expandafter\@slowromancap\romannumeral #1@}
\makeatother

\newcommand{\wt}[1]{\widetilde{#1}}

\newcommand\myeq[1]{\mathrel{\stackrel{\makebox[0pt]{\mbox{\normalfont\tiny #1}}}{=}}}
\newcommand\myge[1]{\mathrel{\stackrel{\makebox[0pt]{\mbox{\normalfont\tiny #1}}}{\ge}}}

\newcommand\mygesim[1]{\mathrel{\stackrel{\makebox[0pt]{\mbox{\normalfont\tiny #1}}}{\gtrsim}}}

\newcommand{\ud}{\,\mathrm{d}}

\newcommand{\prob}{\mathbb{P}}

\newcommand{\Real}{\mathbb{R}}

\newcommand{\ee}{\mathbb{E}}

\newcommand{\vect}[1]{\vec{#1}}

\newcommand{\eps}{\epsilon}

\newcommand{\abs}[1]{\lvert#1\rvert}
\newcommand{\absbig}[1]{\big\lvert#1\big\rvert}
\newcommand{\Abs}[1]{\Big\lvert#1\Big\rvert}

\newcommand{\norm}[1]{\lVert#1\rVert}
\newcommand{\normbig}[1]{\bigl\lVert#1\bigr\rVert}

\newcommand{\inner}[2]{\langle#1, #2\rangle}

\def\bigl{\mathopen\big}

\def\bigr{\mathclose\big}

\newcommand{\unit}{\mathbb{I}}

\newcommand{\detalg}{\mathcal{A}^{\text{Det}}}
\newcommand{\randalg}{\mathcal{A}^{\text{Rand}}}

\newcommand{\emm}{{Euler-Maruyama}}

\newcommand{\probsp}{\mathbb{M}}

\newcommand{\alphaRD}{\eta}
\newcommand{\bmpair}{\wt{W}^{(0,2)}}
\newcommand{\phimap}[2]{\psi^{(#1)}_{#2}}

\newcommand{\revise}[1]{{#1}}
\newcommand{\reviseT}{{T}}
\newcommand{\stepn}{{N_s}}
\newcommand{\timermm}{s}

\begin{document}
 \title{Complexity of randomized algorithms for underdamped Langevin dynamics
 }

\author{Yu Cao\thanks{Department of Mathematics, Duke University, Box 90320, Durham NC 27708, USA
  ({yucao@math.duke.edu}). This author is currently affiliated with Courant Institute of Mathematical Sciences, New York University, USA ({yu.cao@nyu.edu}).}
\and Jianfeng Lu\thanks{Department of Mathematics, Department of Physics, and Department of Chemistry, Duke University, Box 90320, Durham NC 27708, USA ({jianfeng@math.duke.edu}).}
\and Lihan Wang\thanks{Department of Mathematics, Duke University, Box 90320, Durham NC 27708, USA
  ({lihan@math.duke.edu}).}
}

 \pagestyle{myheadings} \markboth{Randomized algorithms for underdamped Langevin}{Y. Cao, J. Lu, and L. Wang}
 \maketitle

\begin{abstract}
We establish an information complexity lower bound of randomized algorithms for simulating underdamped Langevin dynamics. More specifically, we prove that the worst $L^2$ strong error is of order $\Omega(\sqrt{d}\, N^{-3/2})$, for solving a family of $d$-dimensional underdamped Langevin dynamics, by any randomized algorithm with only $N$ queries to $\nabla U$, the driving Brownian motion and its weighted integration, respectively.
The lower bound we establish matches the upper bound for the  randomized midpoint method recently proposed by Shen and Lee [NIPS 2019], in terms of both parameters $N$ and $d$.
\end{abstract}



\section{Introduction}

The underdamped Langevin dynamics have been widely used to sample high-dimensional probability distributions \cite{zwanzig1973nonlinear, mori1965continued,pavliotis_stochastic_2014}, as it could provide a faster convergence rate compared to the overdamped Langevin dynamics. The analysis of the sampling algorithms based on underdamped Langevin dynamics consists of two key aspects:
\begin{enumerate}[(i)]
\item the mixing time of continuous-time underdamped Langevin dynamics;
\item the time-discretization error for numerically integrating underdamped Langevin dynamics.
\end{enumerate}
The first question has been widely studied for various metrics of convergence, see, \eg{}, \cite{pavliotis_stochastic_2014, dolbeault_hypocoercivity_2015, dalalyan_sampling_2020, roussel_spectral_2018, cao_explicit_2019}.

Our focus in this work is the performance of discretization algorithms for underdamped Langevin dynamics. This has also been quite extensively studied, in terms of both asymptotic analysis \cite{pastor_analysis_1988,mishra1996notion,wang2003analysis, gronbech2013simple} and non-asymptotic analysis  \cite{cheng2018sharp,cheng_underdamped_2018,dalalyan_sampling_2020,ma2019there,shen_randomized_2019}. The algorithm with the best rate up to date was proposed by Shen and Lee \cite{shen_randomized_2019}. Their randomized midpoint method \revise{(RMM)} for underdamped Langevin dynamics has a strong $L^2$ error $\mathcal{O}\big(\sqrt{d} N^{-3/2}\big)$ using only
$N$ gradient queries, where $d$ is the dimension.
On the other hand, it is not clear yet from the literature what error rate an optimal algorithm can achieve. In other words, what the intrinsic difficulty of numerical integration of underdamped Langevin dynamics is. This paper provides an answer in the framework of \emph{information-based complexity} (IBC). In particular, we show that the randomized midpoint method is order optimal with respect to both $d$ and $N$.

Information-based complexity  \cite{novak_deterministic_1988,traub_brief_2009}, which is closely related to the notion of the information-theoretic lower bound,  studies the intrinsic complexity of a family of computational problems, based on the type of queries that one has, rather than focusing on a particular algorithm for the task. Intuitively, the algorithmic performance would depend on the information one could acquire (for example, the gradients of the potential function $U$ for Langevin dynamics). IBC aims to establish a lower bound of the accuracy of a family of algorithms, provided the amount and the type of information.

In this work, we adopt the framework of IBC to study randomized algorithms for approximating the strong solution of the underdamped Langevin dynamics with gradient queries to strongly convex potentials and also the driving Brownian motion.

\subsection{Underdamped Langevin dynamics.}
\label{subsec::ULD}
We consider the following underdamped Langevin dynamics $(X_t, V_t) \in \Real^d\times \Real^d$ (we adopt the parameter scaling used in  \cite{cheng_underdamped_2018, shen_randomized_2019}, which is slightly different from the usual physical model of underdamped Langevin dynamics)
\begin{align}
\label{eqn::underdamped}
\begin{split}
\ud X_t &= V_t\ud t,\\
\ud V_t &= - 2 V_t\ud t - \frac{1}{L} \nabla U(X_t)\ud t + \frac{2}{\sqrt{L}}  \ud W_t,
\end{split}
\end{align}
on the time interval $[0, T]$, with the fixed initial condition \revise{$X_0 = x^{\star}$} and $V_0 = 0$, where
\revise{$x^{\star} \equiv x^{\star}(U)\in \Real^d$} is a local minimum of the potential function $U$; \revise{$W_t$ is the $d$-dimensional standard Brownian motion};
the parameter $L>0$ has a physical meaning as the mass of the particle.
The unique stationary distribution of \eqref{eqn::underdamped} is $\rho_{\infty}(x, v) \propto \exp\bigl(-U(x) - \frac{\abs{v}^2}{2/L}\bigr)$. As time $t\to \infty$, the distribution of $(X_t,~V_t)$ converges to the equilibrium exponentially fast under mild conditions; see, \eg{}, \cite{mattingly_ergodicity_2002,villani_hypocoercivity_2009, pavliotis_stochastic_2014, dolbeault_hypocoercivity_2015, roussel_spectral_2018,eberle_couplings_2019, cao_explicit_2019}. Generalization of our main result (\thmref{thm::underdamped_complexity_finite_nablaU}) below to underdamped Langevin dynamics with general friction coefficient is straightforward, and we will not pursue such generality herein for simplicity.

\smallskip
\begin{assumption}
In this work, we shall only consider strongly convex $U$ with Lipschitz gradient, \ie, we consider the following family of potential functions,
\begin{align}
\label{eqn::problem}
\begin{aligned}
\mathcal{F} &\equiv \mathcal{F}(d,\ell,L) :=
\Bigl\{ U\in C^2(\Real^d)\ \Big\rvert\  \ell \unit_d \le \nabla^2 U(x) \le L \unit_d,\ \forall x\in \Real^d
\Bigr\}
\end{aligned}
\end{align}
for fixed parameters $0 < \ell < L < \infty$. Under the strong convexity assumption, we know that $x^{\star}$ is uniquely determined by $U$.
The condition number $\varkappa$ is defined as $\varkappa := L/\ell$.
\end{assumption}
\smallskip

\subsection{Main results.}
\revise{Denote the probability space for underdamped Langevin dynamics \eqref{eqn::underdamped} as $(\probsp, \Sigma, \prob)$.}
In the context of Langevin sampling, we assume that, besides the Brownian motion, the query to a weighted Brownian motion $\wt{W}_t^{(\theta)}$ is also admissible, where
\revise{
\begin{align}
\label{eqn::weighted_bm}
\wt{W}_t^{(\theta)} (\omega) := \int_{0}^{t} e^{\theta s}\ud W_s(\omega), \qquad \forall\ \omega\in \probsp.
\end{align}
When $\theta = 0$, $\wt{W}_t^{(\theta=0)}(\omega)\equiv W_t(\omega)$. 
In general, we define a correlated Gaussian process $\wt{W}_t^{\vect{\theta}}$ for $\vec{\theta} = (\theta_1, \theta_2, \cdots, \theta_{\mathsf{J}})$:
\begin{align*}
\wt{W}_t^{\vect{\theta}}(\omega) := \begin{pmatrix} \wt{W}^{(\theta_1)}_t(\omega) & \cdots & \wt{W}^{(\theta_j)}_t(\omega) & \cdots & \wt{W}^{(\theta_\mathsf{J})}_t(\omega) \end{pmatrix},
\end{align*}
as a short-hand notation. In particular, we shall use $\bmpair_t \equiv (W_t, \wt{W}_t^{(2)})$ frequently below, as in our main theorem.}
The reason for such an assumption is that in the context of Langevin sampling, generating correlated Gaussian random vectors is not computationally expensive, whereas computing $\nabla U$ is usually the computational bottleneck.

\revise{More formally, we assume that there is an oracle query $\Upsilon_{U, \omega}: \Real^d\times [0,T] \rightarrow
\Real^d \times \Real^{2d}$, defined as 
\begin{align*}
    \Upsilon_{U, \omega}(x, t) := \Big(\nabla U(x), \bmpair_t(\omega)\Big),\ x\in \Real^d,\ t\in [0,T],
\end{align*}
for any $U\in \mathcal{F}$ and $\omega\in\probsp$, and the set of \emph{admissible information operations} is defined by
\begin{align}
\label{eqn::acc_info}
    \Lambda := \Big\{ \Upsilon_{U,\omega} (x,t) \big\rvert\ x\in \Real^d,\ t\in [0,T] \Big\}.
\end{align}}

Our main result is the following information-based complexity bound for solving the underdamped Langevin dynamics with $U\in \mathcal{F}$.

\smallskip
\begin{theorem}[Information-based complexity with queries to $\nabla U$ and weighted Brownian motions]
\label{thm::underdamped_complexity_finite_nablaU}
Consider the complexity problem $\mathcal{F}$ \eqref{eqn::problem} \revise{with $\Lambda$ as the set of admissible information operations.}
Then whenever $N \ge N_0$ for some integer $N_0$ (independent of the dimension $d$),
\begin{align}
\label{eqn::main_theorem}
C_{\text{low}} \sqrt{d} N^{-3/2}  \lesssim \inf_{A\in \randalg_{N}} e_{\mathcal{F}, \Lambda} (A) \lesssim C_{\text{up}} \sqrt{d} N^{-3/2},
\end{align}
	where the prefactor $C_{\text{up}} = \sqrt{\frac{T^3}{\ell} + \frac{T^4}{L}}$, and
$C_{\text{low}} \equiv C_{\text{low}}(\ell, L, T)$ can be chosen as
	\begin{align}
	\label{eqn::c_low}
	    C_{\text{low}} = \sup_{\substack{C_x > 0,\ C_v > 0,\\ \ell < u < u_R \le L}} \sqrt{\mathsf{P}(C_x,  C_v, u, T)} C_x^2 \min\{ u-\ell, u_R - u\} \overline{C}(C_x, C_v, u_R, L, T),
	\end{align}
	where $\mathsf{P}(C_x, C_v, u, T)$ for $\ell \le u \le L$ is defined below in \eqref{eqn::P_lower_prob}, and $\overline{C}(C_x, C_v, u_R, L, T)$ for $\ell \le u_R \le L$ is defined below in \eqref{eqn::c_bar}.
\end{theorem}
\smallskip

As a remark, the choice of $\bmpair_t$ \revise{inside the oracle query $\Upsilon_{U,\omega}$ as well as the set of information operations $\Lambda$} (in particular, \revise{the choice of }the weighted Brownian motion $\wt{W}_t^{(2)}$) comes from the friction coefficient in the underdamped Langevin dynamics \eqref{eqn::underdamped}.
In the above, $e_{\mathcal{F}, \Lambda}(A)$ is the worst $L^2$ strong error for any algorithm $A$, defined later in \eqref{eqn::alg_error}; the notation $\randalg_{N}$ means
the set of randomized algorithms that use \revise{$N$ information operations $\{\Upsilon_{U,\omega}(Y_j, t_j)\}_{j=1,2,  \cdots, N}$ (namely, $N$ queries of $\nabla U$ and $N$ queries of $\bmpair_t$)}. The notion of randomized algorithms using only $N$ queries will be elaborated further  in \secref{subsec::ibc} below.
The proof of \thmref{thm::underdamped_complexity_finite_nablaU} will be given in \secref{sec::underdamped}.
The proof of the lower bound estimate is based on a novel non-asymptotic perturbation result with respect to the potential $U$ (see \propref{prop::diff_x_lower}).

\smallskip
\begin{remark}\
\begin{enumerate}[(i)]

\item As a corollary to the theorem, the  randomized midpoint method (see \eqref{eqn::randomized_alg} below for the algorithm) is order optimal.
\item
The fundamental challenge for the computational problem $(\mathcal{F}, \Lambda)$ comes from the insufficient information of $\nabla U$, instead of the path irregularity of the random process $\bmpair_t$. Therefore, $\bmpair_t$ does not play as a bottleneck in proving the lower bound estimate in \thmref{thm::underdamped_complexity_finite_nablaU}. However, if we replace $\bmpair_t$ in the \revise{admissible information operations} $\Lambda$ \eqref{eqn::acc_info} by $W_t$ only (the Brownian motion itself), then the complexity lower bound becomes $\Omega(\sqrt{d} N^{-1})$, as the irregularity of the Brownian motion $W_t$ becomes the complexity bottleneck. This follows from the classical result by Clark and Cameron~\cite{clark_maximum_1980} (see also the literature review in the next subsection).

\item If we replace $\mathcal{F}$ in \eqref{eqn::problem} by the following larger set of potentials
\begin{align*}
\Big\{ U \in C^1(\Real^d)\ \Big\rvert\ \ell \le \frac{\norm{\nabla U(x) - \nabla U(y)}_2}{\norm{x-y}_2}\le L,\ \forall x\neq y\in \Real^d \Big\},
\end{align*}
the above lower bound in \eqref{eqn::main_theorem} also holds, by the definition of $e_{\mathcal{F}, \Lambda}$ \eqref{eqn::alg_error}.

\item In \eqref{eqn::c_low}, the scaling of $C_{\text{low}}$ with respect to the time $T$ and the condition number $\varkappa = L/\ell$ is complicated. Providing a tight estimate of $C_{\text{low}}$ appears to be rather challenging, and we shall leave it to future works.

\item Though we focus on underdamped Langevin dynamics in this paper, the ideas and techniques herein are applicable to a wide range of SDEs.
\end{enumerate}
\end{remark}

\subsection{Literature review.}
Below is a brief literature review for underdamped Langevin algorithms and also the information-based complexity for differential equations.

\subsubsection*{Underdamped Langevin algorithms}

The \emm{} method for the underdamped Langevin dynamics, which replaces $\nabla U(X_t)$ with $\nabla U(X_{kh})$ in \eqref{eqn::underdamped} and solves the modified equation for a short time $h$ at the $k^{\text{th}}$ time step, is the most widely studied algorithm.
To the best of our knowledge, it was first proposed and studied by Ermak and Buckholz  \cite{ermak_numerical_1980}.
Cheng, Chatterji, Bartlett, and Jordan \cite{cheng_underdamped_2018} and later Dalalyan and Riou-Durand \cite{dalalyan_sampling_2020} proved that the strong error of the \emm{} algorithm is $\mathcal{O}(\sqrt{d}N^{-1})$. Results were generalized to non-convex potentials in \cite{cheng2018sharp} with the same rate in $d$ and $N$ but with much worse prefactors. The sampling error of the \emm{} algorithm in Kullback–Leibler divergence was studied in \cite{ma2019there}. Recently, Shen and Lee \cite{shen_randomized_2019} proposed the randomized midpoint method, which reduces the error to $\mathcal{O}(\sqrt{d}N^{-{3}/{2}})$. There have been other algorithms, for example, the BBK scheme \cite{brunger1984stochastic,pastor_analysis_1988}, the Verlet-type scheme proposed in \cite{gronbech2013simple}, and the Leimkuhler-Matthews scheme \cite{leimkuhler2013rational}, but non-asymptotic analysis has not been studied for these yet.

\subsubsection*{Information-based complexity for differential equations}
We now provide a concise review of related works on the information-based complexity for both ordinary and stochastic differential equations.

Information-based complexity analysis for ODEs was initially studied by Kacewicz \cite{kacewicz_optimal_1982,kacewicz_optimality_1983,kacewicz_how_1984,kacewicz_optimal_1987} for deterministic algorithms, and later by Kacewicz \cite{kacewicz_randomized_2004,kacewicz_almost_2006}, Heinrich and Milla \cite{heinrich_randomized_2008}, and Daun \cite{daun_randomized_2011} for randomized algorithms.
For some particular classes of ODE systems, the matching complexity bounds and order optimal algorithms are well known, for both deterministic and randomized algorithms.
A notable observation is that compared to deterministic algorithms, randomized algorithms may achieve order $1/2$ speed-up, which relates to the universal convergence rate of Monte Carlo methods.
This phenomenon also occurs for solving SDE in our case: compared to the \emm{} method, the randomized midpoint method \cite{shen_randomized_2019} achieves order $1/2$ improvement.
On the other hand, we shall comment that it is still open whether such an improvement is non-trivial in our problem by employing randomized algorithms, namely, whether there exists a strong order $3/2$ deterministic algorithm with only gradient queries.

As for the information-based complexity result for SDEs, the lack of full information about both drift and diffusion terms might contribute to the overall computational complexity. It is a common practice to study the complexity from drift and diffusion separately.
Therefore, most works in the literature focus on the complexity due to the diffusion term, since the complexity of the drift term (with trivial diffusion term) reduces to the ODE problem.
However, in our problem, it appears unlikely to reduce the problem into analyzing drift and diffusion terms separately, because both terms are non-trivial for \revise{our computational problem of simulating underdamped Langevin dynamics}. As far as we know, our problem does not fit into known IBC problem formulations for SDEs.

The study of information-based complexity for the diffusion part dates back at least to the seminal work of Clark and Cameron \cite{clark_maximum_1980}. They considered SDEs with $C^3$ drift term (with bounded derivatives up to the third-order) and constant diffusion term, and they proved that the strong one-point approximation error is asymptotically $\Omega(N^{-1})$, for algorithms with uniform mesh grids, where $N$ is the number of queries to the Brownian motion $W_t$ \cite[{Theorem 1}]{clark_maximum_1980}.
The order optimal algorithm is simply the \emm{} method \cite{platen_introduction_1999,higham_algorithmic_2001}.
Clark and Cameron also showed that for general SDEs with Lipschitz diffusion terms, a certain family of algorithms with uniform mesh grid will result in the one-point approximation error with order $\Omega(N^{-1/2})$ \cite{clark_maximum_1980}.
The minimum one-point approximation error for scalar SDEs could be found in \cite{muller-gronbach_optimal_2004}.
Apart from the one-point approximation error, the error for trajectories is also considered, \ie{}, the global approximation error. A series of works of Hofmann, Müller-Gronbach, and Ritter, addressed this problem for $L^{2}$ error \cite{hofmann_optimal_2000,hofmann_optimal_2001,hofmann_global_2004}, and for
$L^{\infty}$ error \cite{hofmann_step_2000, muller-gronbach_optimal_2002}.
We refer readers to a survey paper by Müller-Gronbach and Ritter \cite{muller-gronbach_minimal_2008} for more details.
Around the last decade, Przybyłowicz \cite{przybylowicz_adaptive_2010,przybylowicz_minimal_2015,przybylowicz_optimal_2015}, Przybyłowicz and Morkisz \cite{przybylowicz_strong_2014} studied the time-irregular SDEs.
More recently, Hefter, Herzwurm, and Müller-Gronbach provided probabilistic lower bound estimates \cite{hefter_lower_2019}.

Finally, we also point out other related works studying differential equations with inexact information, see, \eg{}, \cite{kacewicz_optimal_2016,bochacik_randomized_2020} for ODEs and \cite{morkisz_optimal_2017,morkisz_randomized_2021} for SDEs.

\subsection*{Notation} The Lebesgue measure on $\Real^d$ (for any dimension $d$) is denoted by $\mu$. The probability space for \eqref{eqn::underdamped} is $(\probsp, \Sigma, \prob)$ and \revise{the probability space of the source of randomness for randomized algorithms is denoted by $(\wt{\probsp}, \wt{\Sigma}, \wt{\prob})$.
Suppose $\mathcal{X}, \mathcal{Y}$ are two arbitrary spaces, and let $\mathscr{X}: \probsp\times \mathcal{X}\rightarrow \mathcal{Y}$ be a function such that $\mathscr{X}(\cdot, \iota)$ is a random variable on $(\probsp, \Sigma, \prob)$ for any fixed $\iota\in \mathcal{X}$. Then,  
\begin{align*}
\ee_{{\omega}}[\mathscr{X}(\omega, \iota)] \equiv \ee_{{\omega}\sim \prob} [\mathscr{X}(\omega, \iota)] := \int \mathscr{X}(\omega, \iota)\ \ud\prob(\omega).
\end{align*}  
Other notations like $\ee_{\wt{\omega}} \equiv \ee_{\wt{\omega}\sim \wt{\prob}}$ and $\ee_{(\omega, \wt{\omega})} \equiv \ee_{(\omega, \wt{\omega})\sim \prob\times \wt{\prob}}$ are similarly defined.
}

\revise{For two arbitrary real-valued functions $f$ and $g$,}
when $c f\le g$ for some universal constant $c$,
we denote this by $f\lesssim g$ as a simplification; the notation $\gtrsim$ is similarly defined.
The notations $\Omega$ and $\mathcal{O}$ follow the convention in complexity analysis, \ie{}, $f = \Omega(g)$ means $f \gtrsim g$, and $f = \mathcal{O}(g)$ means $f\lesssim g$.
The binary relations $\prec$ and $\succ$ are partial order relations on the Boolean lattice $\{0,1\}^N$.

\section{Preliminaries}

In this section, we will illustrate the setup of the IBC problem under consideration,
the integral form of \eqref{eqn::underdamped}, the randomized midpoint method \cite{shen_randomized_2019}, and the exact solution of \eqref{eqn::underdamped} for $1D$ quadratic potentials.

\subsection{IBC problem setup.}
\label{subsec::ibc}
We have explained the family of computational problem, characterized by the family $\mathcal{F}$ \eqref{eqn::problem} in the introduction. Next, we shall explain more about the \revise{admissible information operations} and the family of randomized algorithms. We refer readers to \eg{}, \cite[Sec. 2]{muller-gronbach_minimal_2008} for a more abstract framework.

\revise{
\subsubsection*{Solution map}
For the computational problem under consideration, we want to approximate the
following solution map, for a given $T> 0$:
\begin{align*}
    X_T: (U, \omega)\in (\mathcal{F}\times \probsp) \mapsto X_T(U,\omega)\in \Real^d,
\end{align*}
where $X_T(U, \omega)$ denotes the strong solution  of \eqref{eqn::underdamped} at time $T$, for a given potential function $U$.

\subsubsection*{Admissible information operations}

The set of admissible information operations $\Lambda$ \eqref{eqn::acc_info}, recalled here, is  
\begin{align*}
    \Lambda = \Bigl\{ \Upsilon_{U,\omega} (x,t) \;\big\vert\; x\in \Real^d,\ t\in [0,T] \Bigr\}.
\end{align*}
We need to use finite amount of information (operations) to predict the strong solution $X_T(U, \omega)$ in \eqref{eqn::underdamped}, and such a prediction is known as an algorithm.

\subsubsection*{Algorithms}

We shall always choose 
\begin{align}
\label{eqn::Y1_t1}
    Y_1 = x^{\star}(U),\qquad t_1 = 0,
\end{align}
as these are the
initial conditions in the dynamics \eqref{eqn::underdamped}. 
A deterministic algorithm $A$ with $N$ evaluations means there is a sequence of deterministic mappings $\varphi_j: \Real^d\times (\Real^d\times \Real^{2d})^{j}\rightarrow \Real^d\times \Real$ ($1\le j \le N-1$) 
and a mapping $\phi :\Real^d\times (\Real^d\times\Real^{2d})^{N}\rightarrow\Real^d$ such that
\begin{subequations}
\label{eqn::det_alg_maps}
    \begin{align}
    (Y_{j+1}, t_{j+1}) := \varphi_j\big(Y_1, \Phi_j\big) \equiv \varphi_j\bigl(x^{\star}(U), \Phi_j\big),\qquad \text{ for }\ 1\le j \le
    N-1;\\
    \Phi_j := \big(\Upsilon_{U,\omega}(Y_1, t_1),
        \Upsilon_{U,\omega}(Y_2, t_2), \cdots,
    \Upsilon_{U,\omega}(Y_{j},t_j)\big),\qquad
        \text{ for } 1\le j \le N;
    \label{eqn::info_j}
    \end{align}
\end{subequations}
and the algorithm $A$ is given by 
\begin{align}
\label{eqn::alg_form}
\begin{aligned}
A(U, \omega) &= \phi(Y_1, \Phi_{N}) \\
&\equiv \phi\big(x^{\star}(U), \Upsilon_{U,\omega}(Y_1, t_1), \Upsilon_{U,\omega}(Y_2,
    t_2), \cdots, \Upsilon_{U, \omega}(Y_{N}, t_{N})\big).
\end{aligned}
\end{align}
Ideally, we hope that $X_T(U,\omega)\approx A(U,\omega)$ under certain norms.
The family of all such algorithms is denoted by $\detalg_{N}$.

\smallskip

The family of randomized algorithms with only $N$ queries is denoted by $\randalg_{N}$. 
Consider any other probability space $(\wt{\probsp}, \wt{\Sigma}, \wt{\prob})$
as the source of randomness. 
We define randomized algorithms $A( U,\omega, \wt{\omega})$ as
mappings from the space $\mathcal{F} \times
\probsp\times \wt{\probsp}$ to the vector space $\Real^d$:
\begin{align*}
 (U, \omega, \wt{\omega}) \in \mathcal{F}\times \probsp \times \wt{\probsp}
\rightarrow A(U,\omega, \wt{\omega}) \in \Real^d,
\end{align*}
such that 
\begin{align*}
 A(\cdot,\cdot, \wt{\omega})\in \detalg_{N},\qquad    \forall\  \wt{\omega}\in \wt{\probsp}.
\end{align*}
Apparently, $\detalg_{N} \subseteq \randalg_{N}$ and moreover, $\detalg_{N_1} \subseteq \detalg_{N_2}$,  $\randalg_{N_1} \subseteq \randalg_{N_2}$ whenever $N_1 \le N_2$.

The error of the randomized algorithm $A(\cdot,\cdot,\wt{\omega})$ is measured in the $L^2$ sense herein:
\begin{align}
\label{eqn::alg_error}
e_{\mathcal{F}, \Lambda}(A) := \sup_{U \in \mathcal{F}} \Bigl(\ee_{(\omega,
\wt{\omega})\sim \prob\times \wt{\prob}}\bigl[ \abs{X_T(U, \omega) - A(U, \omega, \wt{\omega})}^2\bigr]\Bigr)^{1/2}.
\end{align}
}

\smallskip
\begin{remark}
For the underdamped Langevin dynamics used to sample log-concave probability distributions $e^{-U}/\int e^{-U} $, we are interested only in $X_T$, instead of the whole trajectory on the interval $[0,T]$.
\end{remark}

\subsection{The integral form and the randomized midpoint method.}
\label{subsec::numerical_alg}
Several numerical algorithms for underdamped Langevin dynamics \cite{ermak_numerical_1980,cheng_underdamped_2018,dalalyan_sampling_2020,shen_randomized_2019} are based on  its integral form
\begin{align}
\label{eqn::integral_form}
\begin{aligned}
& \begin{aligned}
X_t =
X_0 + \frac{1-e^{- 2 t}}{2} V_0  + & \frac{1}{\sqrt{L}} \int_{0}^{t} (1-e^{2 (s-t)}) \ud W_s
-\frac{1}{2L} \int_{0}^{t} (1-e^{2(s-t)}) \nabla U(X_s)\ud s,\\
\end{aligned}\\
& \begin{aligned}
V_t & = e^{-2 t} V_0 + \frac{2}{\sqrt{L}} \int_{0}^{t} e^{2(s-t)} \ud W_s  -\frac{1}{L} \int_{0}^{t} e^{2 (s-t)}\nabla U(X_s)\ud s.
\end{aligned}
\end{aligned}
\end{align}

In \cite{shen_randomized_2019}, Shen and Lee considered a randomized midpoint method for simulating underdamped Langevin dynamics, in the context of sampling log-concave distributions.
\revise{Let $\stepn$ be the total number of steps and $h :=T/\stepn$ is the time step size. Given $\hat{X}_k$ and $\hat{V}_k$ at time $\timermm_k := k h$ for integer $0\le k \le \stepn-1$,} the randomized midpoint method approximates $\hat{X}_{k+1}$ and $\hat{V}_{k+1}$ in the following way:\revise{
\begin{align}
\label{eqn::randomized_alg}
\begin{aligned}
\hat{X}_{k+1} & = \hat{X}_{k} + \frac{1-e^{-2 h}}{2} \hat{V}_k + \frac{1}{\sqrt{L}}\int_{0}^{h} (1-e^{2(s-h)})\ud W_{\timermm_k+s} \\
& \qquad - \frac{1}{2L} h (1-e^{2(\alphaRD_k h - h)})\nabla U(\hat{X}_{k+1/2}),\\
\hat{V}_{k+1} & = e^{-2 h} \hat{V}_k + \frac{2}{\sqrt{L}}\int_{0}^{h} e^{2(s - h)}\ud W_{\timermm_k+s}  \\
&\qquad - \frac{1}{L} h e^{2(\alphaRD_k h - h)}\nabla U(\hat{X}_{k+1/2}),\\
\hat{X}_{k+1/2} & = \hat{X}_k +  \frac{1-e^{-2 \alphaRD_k h}}{2} \hat{V}_k + \frac{1}{\sqrt{L}}\int_{0}^{\alphaRD_k h} (1-e^{2(s-\alphaRD_k h)})\ud W_{\timermm_k+s} \\
&\qquad -\frac{1}{2L} \big(\int_{0}^{\alphaRD_k h} 1-e^{2(s-\alphaRD_k h)}\ud s\bigr) \nabla U(\hat{X}_k),
\end{aligned}
\end{align}
where \emph{i.i.d.} random variables $\bigl\{\alphaRD_k\bigr\}_{0\le k \le \stepn-1}$ are uniformly distributed on the interval $[0,1]$, independent of $W_{t}$; the initial condition is $\hat{X}_0 = x^{\star}(U)$ and $\hat{V}_0 = 0$.}

In the above, $\hat{X}_{k+1}$ and $\hat{V}_{k+1}$ are obtained by approximating the integral with respect to $\nabla U(X_s)$ in \eqref{eqn::integral_form} by its evaluation at a single point $\hat{X}_{k+1/2}$; $\hat{X}_{k+1/2}$ is obtained using the \emm{} scheme for the integral form \eqref{eqn::integral_form} at the random time \revise{$\timermm_k + \alphaRD_k h$}. 

\smallskip
\revise{
\begin{remark}
The RMM with $\stepn$ steps is a randomized algorithm with $N = 2 \stepn$ evaluations in the form of \eqref{eqn::det_alg_maps} and \eqref{eqn::alg_form}; please refer to Appendix \ref{app::rmm-randomized} for more details.
\end{remark}}
\smallskip

Many previous works have proposed and analyzed various randomized algorithms to solve differential equations.
For ODEs, the analogy of the randomized midpoint method can be found in, \eg{}, \cite{heinrich_randomized_2008,daun_randomized_2011,bochacik_randomized_2020}.
More randomized ODE solvers could be found in \eg{}, \cite{jentzen_random_2009,stengle_numerical_1990,stengle_error_1995}.
As for SDEs, the randomized Euler's method, which is in a very similar spirit as the randomized midpoint method, has been studied in, \eg{}, \cite{przybylowicz_strong_2014,przybylowicz_minimal_2015,przybylowicz_optimal_2015,morkisz_optimal_2017}.
A randomized Milstein method was studied by Kruse and Wu
\cite{kruse_randomized_2018} for non-differentiable drift functions; a randomized derivative-free Milstein method was studied by Morkisz and Przybyłowicz \cite{morkisz_randomized_2021} for scalar SDEs with inexact information.

\smallskip
\begin{remark}\normalfont
While we focus on the strong error in this paper, we would like to comment that for underdamped Langevin dynamics, the randomized midpoint method defined in \eqref{eqn::randomized_alg} has a weak error of $\mathcal{O}(h^3)$. As the full analysis is tedious and it is not the main focus of our paper, we only include a heuristic argument here. \revise{We shall consider a one-step error only and let $\eta \equiv \eta_0$ herein for simplicity.}

Note that for a fixed realization of $W_t$, we can compare
\begin{align*}
   \hat{X}_1-X_h & = -\dfrac{1}{2L}\Bigl(h(1-e^{-2(h-\alphaRD h)})\nabla U(\hat{X}_\frac{1}{2})-\int_0^h (1-e^{-2(h-s)})\nabla U(X_s)\ud s \Bigr) \\ & = -\dfrac{1}{2L}\Bigl(h(1-e^{-2(h-\alphaRD h)})\nabla U(\hat{X}_\frac{1}{2})-h(1-e^{-2(h-\alphaRD h)})\nabla U(X_{\alphaRD h})\\ & \qquad +h(1-e^{-2(h-\alphaRD h)})\nabla U(X_{\alphaRD h}) -\int_0^h (1-e^{-2(h-s)})\nabla U(X_s)\ud s \Bigr),
\end{align*}
\begin{align*}
    \hat{V}_1-V_h & = -\dfrac{1}{L}\Bigl(he^{-2(h-\alphaRD h)}\nabla U(\hat{X}_\frac{1}{2})-\int_0^h e^{-2(h-s)}\nabla U(X_s)\ud s \Bigr) \\
    & = -\dfrac{1}{L}\Bigl(he^{-2(h-\alphaRD h)}\nabla U(\hat{X}_\frac{1}{2})-he^{-2(h-\alphaRD h)}\nabla U(X_{\alphaRD h}) \Bigr)\\ & \qquad -\frac{1}{L}\Bigl(he^{-2(h-\alphaRD h)}\nabla U(X_{\alphaRD h}) -\int_0^h e^{-2(h-s)}\nabla U(X_s)\ud s \Bigr) \\ &=: I+II.
\end{align*}
We observe that $\hat{X}_1-X_h$ is of higher order than $\hat{V}_1-V_h$. Therefore, keeping only the error from $V$ and by the Taylor expansion around $(X_h,V_h)$, we can estimate, for a smooth enough test function $f$, that
\begin{align*}
   \absbig{\ee \bigl[ f(\hat{X}_1,\hat{V}_1)-f(X_h,V_h)\bigr] }  &\lesssim \absbig{\ee \bigl[\nabla_v f(X_h,V_h)\cdot (I+II) \bigr] } \\ & \qquad  + \absbig{\ee \bigl[(I+II) \cdot \nabla^2_v f(X_h,V_h) (I+II) \bigr] }.
\end{align*} 
\revise{Now let $\ee_\eta$ denote the expectation with respect to $\eta$ only. Since $\eta$ is uniformly distributed on $[0,1]$, we have \begin{equation*}
    \ee_\eta \bigl[ he^{-2(h-\alphaRD h)}\nabla U(X_{\alphaRD h})\bigr] = h\int_0^1 e^{-2(h-\alphaRD h)}\nabla U(X_{\alphaRD h}) \ud \alphaRD = \int_0^h e^{-2(h-s)} \nabla U(X_s)\ud s.
\end{equation*} This shows $\ee_\eta \bigl[II\bigr]=0$. Since $(X_h,V_h)$ is independent of $\alphaRD$, we have \begin{equation*}
   \ee \bigl[\nabla_v f(X_h,V_h)\cdot II \bigr] = \ee \Big[\nabla_v f(X_h,V_h) \cdot \ee_\eta \bigl[ II\bigr]\Big] =0.
\end{equation*}Therefore, }
\begin{align*}
     \absbig{\ee \bigl[ f(\hat{X}_1,\hat{V}_1)-f(X_h,V_h)\bigr] } & \lesssim \dfrac{h}{L}\Abs{\ee \bigl[\nabla_v f(X_h,V_h)\cdot \bigl(\nabla U(\hat{X}_\frac{1}{2})-\nabla U(X_{\alphaRD h})\bigr)\bigr] } \\ & \qquad + \ee \bigl[\norm{\nabla^2_v f(X_h,V_h)}\cdot \abs{ \hat{V}_1-V_h}^2 \bigr] .
\end{align*}
Finally, assuming that all derivatives of $f$ are bounded, we use \cite[Lemma 9]{shen_randomized_2019} for the first term below, and \cite[Lemma 2]{shen_randomized_2019} for the second term below,
\begin{align*}
 \absbig{\ee \bigl(&f(\hat{X}_1,\hat{V}_1)-f(X_h,V_h)\bigr)} \lesssim \dfrac{h}{L}\Big(\ee \abs{\nabla U(\hat{X}_\frac{1}{2})-\nabla U(X_{\alphaRD h})}^2\Big)^\frac{1}{2}+\ee \abs{ \hat{V}_1-V_h}^2 =\mathcal{O}(h^4).
 \end{align*}
This local truncation error gives $\mathcal{O}(h^3)$ weak error by Gr{\" o}nwall's inequality as usual.
\end{remark}

\subsection{Exact solution for quadratic potentials in 1D.}
Our estimate of the prefactor $C_{\text{low}}$ in
\thmref{thm::underdamped_complexity_finite_nablaU} relies on the
\revise{behavior of $X_t(U_u, \omega)$ and $V_t(U_u, \omega)$}, where the potential has the quadratic form $U_u(x) := u x^2/2$ (thus \revise{$x^{\star}(U_u) = 0$} for this case).
Therefore, in this subsection, we shall first review the exact solution of the underdamped Langevin dynamics \revise{with quadratic potentials and then define a quantity in \eqref{eqn::P_lower_prob} to be used in the next section.}

It is easy to rewrite \eqref{eqn::underdamped} as
\begin{align*}
    \ud \begin{bmatrix}
    X_t\\ V_t
    \end{bmatrix} = H \begin{bmatrix}
    X_t\\ V_t
    \end{bmatrix} \ud t + \begin{bmatrix}
    0 \\
    \frac{2}{\sqrt{L}}
    \end{bmatrix} \ud W_t,\qquad H = \begin{bmatrix}
    0 & 1 \\
    -\frac{u}{L} & -2 \\
    \end{bmatrix}.
\end{align*}
Then its integral form can be immediately obtained as follows
\begin{align*}
    \begin{bmatrix}
    X_t\\ V_t
    \end{bmatrix} = \int_{0}^{t} e^{H(t-s)} \begin{bmatrix}0\\ \frac{2}{\sqrt{L}} \end{bmatrix}\ud W_s.
\end{align*}
The matrix exponential of $H$ can be explicitly computed \revise{for $u<L$:}
\revise{
\begin{align*}
    e^{H t} = \frac{1}{2\sqrt{1-u/L}}\begin{bmatrix} 
        \Big(\lambda_{+} e^{-t \lambda_{-}} - \lambda_{-} e^{-t
        \lambda_{+}}\Big) &  
    \Big(e^{-t\lambda_{-}} - e^{-t \lambda_{+}}\Big)\\
        \frac{u}{L} \Big(e^{-t\lambda_{+}} - e^{-t\lambda_{-}}\Big) &
    \Big(\lambda_{+} e^{-t\lambda_{+}} - \lambda_{-} e^{-t\lambda_{-}}\Big)\\
\end{bmatrix},
\end{align*}
where
\begin{align}
\label{eqn::lambda}
   \lambda_{\pm} =  1 \pm \sqrt{1-u/L}.
\end{align}
This immediately leads to the following result.
}
\smallskip
\begin{lemma}[Exact solution]
When $d = 1$ and $U_u(x) = u x^2/2$, we have \revise{for $u < L$}
\begin{align}
\label{eqn::exact_soln}
\begin{aligned}
\revise{X_t(U_u, \cdot)} & = \frac{1}{\sqrt{L - u}}  \int_{0}^{t} \big(e^{-(t-s)\lambda_{-}} - e^{-(t-s)\lambda_{+}}\big)\ud W_s, \\
\revise{V_t(U_u, \cdot)} &= \frac{1}{\sqrt{L - u}} \int_{0}^{t}
\big(-\lambda_{-}e^{-(t-s)\lambda_{-}} + \lambda_{+} e^{-(t-s) \lambda_{+}}
\big)\ud W_s.
\end{aligned}
\end{align}
\end{lemma}
\smallskip

Next, let us introduce the following quantity for $C_x, C_v > 0$, $\ell \le u \le L$,
\begin{align}
\label{eqn::P_lower_prob}
\begin{aligned}
    \mathsf{P}(C_x, C_v, u, T)
    :=\ \prob\Big(\omega: &\sup_{0\le t\le T} X_t(U_{u}, \omega) \ge 2C_x,\\  & \inf_{0\le t\le T}
    X_t(U_{u}, \omega) \le - 2C_x, \sup_{0\le t\le T} \abs{V_t(U_{u}, \omega)} \le C_v/2\Big).
  \end{aligned}
\end{align}
The event under consideration requires $X_t(U_u, \omega)$  \revise{to cross $-2
C_x$ and $2 C_x$}, 
whereas the velocity $V_t(U_u, \omega)$ is uniformly bounded on the whole interval $[0, T]$. Typically, we should expect to choose a small $C_x$ and a large $C_v$ in order to have $\mathsf{P}(C_x, C_v, u, T) = \mathcal{O}(1)$.
Indeed, if $C_x / C_v$ is too large, \revise{the probability $\mathsf{P}(C_x, C_v, u, T)$ equals to $0$}, as stated in the following lemma.

\smallskip
\begin{lemma}
When $ 12 C_x/C_v > T$, then $\mathsf{P}(C_x, C_v, u, T) = 0$.
\end{lemma}
\smallskip

\begin{proof}
For any $\omega$ satisfying the condition in \eqref{eqn::P_lower_prob}, the travel distance of $X_t$ must be at least $6 C_x$ (since it starts from $0$ and has to cross levels $2 C_x$ and $-2C_x$) with velocity at most $C_v/2$, then the total time must be at least $12 C_x/C_v$.
\hfill
\end{proof}

\smallskip
\revise{
However, note that for large enough $C_v$, $\mathsf{P}(C_x, C_v, u, T) > 0$.
\smallskip

\begin{lemma}
\label{lem::non-trivial-P}
For fixed $C_x > 0$, $u \in (\ell, L)$, and $T > 0$, and for sufficiently large $C_v$, we have 
\begin{align*}
\mathsf{P}(C_x, & C_v, u, T) \ge \frac{1}{2} \prob\Big(\omega: \sup_{0\le t\le T} X_t(U_{u}, \omega) \ge 2C_x,\  \inf_{0\le t\le T} X_t(U_{u}, \omega) \le - 2C_x\Big) > 0.
\end{align*}
\end{lemma}
The proof of this lemma is postponed to Appendix \ref{app::supplementary}.
}

 Therefore, the prefactor $C_{\text{low}}$ in \eqref{eqn::c_low} is non-zero. The precise dependence of $C_{\text{low}}$ on parameters $\ell$, $L$ and $T$, however, appears to be a challenging problem and will be left for future investigations.

\section{Proof of \thmref{thm::underdamped_complexity_finite_nablaU}}
\label{sec::underdamped}

The proof of the lower bound estimate relies on the (non-asymptotic) perturbation analysis in \secref{sec::perturbation},
in particular, the lower bound estimate in \propref{prop::diff_x_lower}.
The overall strategy, from the information-based complexity perspective, is similar to the lower bound estimate for randomized algorithms for integration problems, see, \eg{}, \cite{novak_deterministic_1988}. The new ingredient is the perturbation type analysis for the particular problem under consideration.
The proof of the upper bound is known from  \cite{shen_randomized_2019}.
Thus, we shall only provide a sketch of the main steps to prove the upper bound for completeness.

\subsection{Non-asymptotic perturbation analysis with respect to $U$.}
\label{sec::perturbation}
We consider the case $d = 1$, which is assumed throughout this subsection. We postpone proofs for all results in this subsection to \secref{subsec::proof_perturbation_lemmas} for clarity.

Let us consider
\begin{align*}
\revise{U_u(x)} := \frac{u x^2}{2},\ \forall u\in [\ell, L],
\end{align*}
and let us also introduce a set parameterized by $u\in (\ell, L)$ and $\eps>0$
\begin{align}
\label{eqn::Fueps}
\mathcal{F}_{u, \eps} :=  \{U\in \mathcal{F}:\ \norm{\nabla U(x) - {u} x}_{\infty} \le \eps,\ x^{\star}(U) = 0 \}.
\end{align}
\revise{Recall from \secref{subsec::ULD} that $x^{\star}(U)$ is the local minimum of potential function $U\in \mathcal{F}$; thus} for any $U\in \mathcal{F}_{u, \eps}$, we have the initial condition $X_0(U,\omega) = 0$.

First, we show that when the potential function is slightly perturbed away from a quadratic function $U_u$, the strong solutions of $X_t(U, \omega)$ and $V_t(U, \omega)$ are at most perturbed by an order of $\mathcal{O}(\eps)$.

\smallskip
\begin{lemma}[Upper bound]
\label{lem::perturbation}
Consider any $u\in (\ell, L)$. For any $U\in \mathcal{F}_{u, \eps}$, any $\omega\in \probsp$, and $t\in [0,T]$, we have
\begin{align}
\label{eqn::diff_xv_upper}
\begin{aligned}
 \abs{X_t(U, \omega) - X_t(U_u, \omega)} &\le \frac{\eps}{2 L (1-\sqrt{1-\frac{u}{L}}) \sqrt{1-\frac{u}{L}}},\\
 \abs{V_t(U, \omega) - V_t(U_u, \omega)} &\le \frac{\eps}{L\sqrt{1-\frac{u}{L}}}.
 \end{aligned}
\end{align}
\end{lemma}
\smallskip

For any $u\in (\ell, L)$, $C_x > 0$ and $C_v > 0$, let us define
\begin{align*}
\begin{aligned}
\bar{\epsilon} &\equiv \bar{\epsilon}(C_x, C_v, u, L) \\
&:=
\min\Big\{2 L \Big(1-\sqrt{1-\frac{u}{L}}\Big) \sqrt{1-\frac{u}{L}} C_x,\ L\sqrt{1-\frac{u}{L}} \frac{C_v}{2}\Big\},
\end{aligned}
\end{align*}
which is strictly positive. Moreover, for $\eps > 0$,
let us define a set
$\mathcal{E} \equiv \mathcal{E}(C_x, C_v, u, T, \eps)$
\begin{align}
\label{eqn::event_E}
\begin{aligned}
\mathcal{E}:= \Big\{\omega: \sup_{0\le t\le T} X_t(U, \omega) \ge C_x,\  \inf_{0\le t\le T} & X_t(U, \omega)  \le - C_x, \\
& \sup_{0\le t\le T} \abs{V_t(U, \omega)} \le C_v,\  \forall U\in \mathcal{F}_{u, \eps}\Big\}.
\end{aligned}
\end{align}

In the following lemma, we shall bound $\prob(\mathcal{E})$ by $\mathsf{P}(C_x, C_v, u, T)$ from below, for small enough $\eps$.

\smallskip
\begin{lemma}
\label{lem::choose_bounds}
Consider any $u\in (\ell,L)$,
and any $C_x, C_v > 0$.
For arbitrary $0 < \eps < \bar{\eps}$, $\prob(\mathcal{E})$ is uniformly bounded from below by $\mathsf{P}(C_x, C_v, u, T)$ defined in  \eqref{eqn::P_lower_prob}.
\end{lemma}
\smallskip

In general, for two fixed potential functions $U_1 \neq U_2$, the distance $\abs{X_T(U_1, \omega) - X_T(U_2, \omega)}$ highly depends on the realization of the Brownian motion \revise{$W_t(\omega)$}, and it is unlikely to establish a uniform non-trivial lower bound of $\abs{X_T(U_1, \omega) - X_T(U_2, \omega)}$ for arbitrary $\omega$. However, if we restrict the outcome $\omega$ to a \enquote{nice} set, \ie, $\omega \in \mathcal{E}$ defined in \eqref{eqn::event_E}, then we can provide a lower bound estimate of $\abs{X_T(U_1, \omega) - X_T(U_2, \omega)}$ as in the following Proposition.  This lower bound estimate is the key to prove \thmref{thm::underdamped_complexity_finite_nablaU}.

\smallskip
\begin{proposition}
\label{prop::diff_x_lower}
\revise{Consider any $u\in(\ell,L)$, $C_x, C_v, \eps > 0$, and let $\mathcal{E} \equiv \mathcal{E}(C_x, C_v, u, T, \eps)$ be the set defined in \eqref{eqn::event_E}.}
Consider two potential functions $U_1, U_2\in \mathcal{F}_{u, \eps}$.
Assume that
\begin{enumerate}[(i)]

\item the continuous function $\revise{g} := \nabla U_1 - \nabla U_2$ is non-negative on $\Real$;

\item \label{assump::support_I} there exists $\mathcal{I}\subseteq [-C_x/2, C_x/2]$, a finite union of closed bounded intervals, such that
\begin{align}
\label{eqn::assump_g}
g(x) \ge \frac{\epsilon}{2} \unit_{\mathcal{I}}(x),\qquad \forall x\in \Real,
\end{align}
\revise{where $\unit_{\mathcal{I}}$ is the indicator function of set $\mathcal{I}$ on $\Real$.}
\end{enumerate}

Let us introduce $u_R\in [u, L]$ as a constant such that
\begin{align*}
\nabla^2 U_2(x) \le u_R, \qquad \forall x\in \Real.
\end{align*}
Then for any $\omega\in \mathcal{E}$, we have
\begin{align*}
\abs{X_T(U_1, \omega) - X_T(U_{2}, \omega)} \ge \overline{C}\ \eps \mu(\mathcal{I}),
\end{align*}
where $\mu$ is the Lebesgue measure and $\overline{C}\equiv \overline{C}(C_x, C_v, u_R, L, T)$ is given by
\begin{align}
\label{eqn::c_bar}
\overline{C} := \left\{\begin{aligned}
&\frac{e^{(\frac{3C_x}{2C_v} - T) (1-\sqrt{1-\frac{u_R}{L}})} \Big(1 - e^{- \frac{C_x}{C_v}\sqrt{1-\frac{u_R}{L}}}\Big)}{4 L C_v\sqrt{1-\frac{u_R}{L}}},\ &\text{ if } u_R < L; \\
&\frac{C_x e^{\frac{3 C_x}{2 C_v} - T}}{4L C_v^2},\ & \text{ if } u_R = L.
\end{aligned}\right.
\end{align}
\end{proposition}

\subsection{Proof of the lower bound estimate for \thmref{thm::underdamped_complexity_finite_nablaU}.}
\label{subsec::proof_lower_nablaU}
We shall proceed to prove the lower bound estimate, based on the results in \secref{sec::perturbation}.

\medskip
{\noindent \bf Case (\rom{1}): $d = 1$.} We shall first consider the case $d = 1$.

\medskip
{\noindent \emph{Step (1): Setup and notations.}}
Since we only access $\nabla U$ at $N$ points, we will not be able to gain the full information of $\nabla U$ based on the local queries.
In this step, we shall consider a family of $U$ (see \eqref{eqn::U_family} and \eqref{eqn::Gset} below) as small perturbations of a quadratic potential with the mode (Hessian) $u$, and we shall estimate how the deviation of $\nabla U$ contributes to the error $e_{\mathcal{F}, \Lambda}(A)$ \eqref{eqn::alg_error} for any randomized algorithm $A\in\randalg_{N}$.

Without loss of generality, we assume that $N$ is an even integer.
We shall pick $\ell <  u <  u_R \le L$. We could then fix $C_x$ and $C_v$ satisfying $\mathsf{P}(C_x, C_v, u, T) > 0$.
Let us define
\begin{align}
\label{eqn::parameters}
\xi  := \min\{ u-\ell, u_R - u\} > 0, \qquad \epsilon := \frac{C_x \xi}{8 N}.
\end{align}
We will pick $N$ sufficiently large such that $\eps < \bar{\epsilon}(C_x, C_v, u, L)$. Then we could also determine the set $\mathcal{E}$ defined in \eqref{eqn::event_E}.

We divide the interval $I := \big[-\frac{C_x}{2}, \frac{C_x}{2}\big]$ into $2N$ equally spaced sub-intervals. Let $x_j = \frac{C_x}{2N} j$ for integers $-N\le j\le N$; \revise{let $I_j = \big[x_j, x_{j+1}\big)$ for $-N\le j \le N-2$ and let $I_{N-1} = \big[x_{N-1}, x_{N}\big]$.}
Define a non-negative function \revise{$g$ supported on $\big[0, \frac{C_x}{2N} \big]$} by
\begin{equation*}
g(x) :=
\begin{cases}
a x^2, \qquad  & x\in \bigl[0, \frac{C_x}{8N} \bigr];\\
- a (x - \frac{C_x}{4N})^2 + 2a \bigl(\frac{C_x}{8 N}\bigr)^2, \qquad & x \in \bigl[ \frac{C_x}{8N},  \frac{3 C_x}{8N}\bigr];\\
a(x - \frac{C_x}{2N})^2, \qquad & x\in \bigl[\frac{3C_x}{8N}, \frac{C_x}{2N}\bigr]; \\
\revise{0,} \qquad & \revise{x\notin \bigl[0, \frac{C_x}{2N}\bigr]},\\
\end{cases}
\end{equation*}
where $a := \xi \frac{4N}{C_x}$. It is easy to verify that $g\in C^1(\Real)$ and moreover,
\begin{align*}
\norm{g'}_{\infty} &= \xi, \qquad \norm{g}_{\infty} = \eps.
\end{align*}
Moreover, we choose $\mathscr{I} = \big[\frac{C_x}{8N}, \frac{3C_x}{8N}\big]$ (thus, the length $\mu(\mathscr{I}) = \frac{C_x}{4N}$), and we know that $g(x) \ge \frac{\eps}{2} \unit_{\mathscr{I}}(x)$ \revise{for any $x\in \Real$}.

For index ${\beta} = (\beta_{-N}, \beta_{-N+1}, \cdots, \beta_{N-1}) \in \{0,1\}^{2N}$, we define $U_{\beta}$ with $U_{\beta}(0) = 0$ by its derivative
\begin{align}
\label{eqn::U_family}
\nabla U_{\beta}(x) := u x + \sum_{j=-N}^{N-1} \beta_j g(x - x_j) \ge u x + \sum_{j=-N}^{N-1} \beta_j \frac{\eps}{2} \unit_{\mathscr{I}}(x-x_j).
\end{align}
Apparently, $U_{\beta}$ is well-defined, $U_{\beta}\in \mathcal{F}_{u,\eps}$, and $\nabla^2 U_{\beta}(x) \le u_R$, $\forall \beta$, $\forall x\in \Real$.
Define a space
\begin{align}
\label{eqn::Gset}
\mathcal{G} \equiv \mathcal{G}_{u, \eps} := \{U_{\beta}:\ \beta \in \{0,1\}^{2N} \}\subseteq \mathcal{F}_{u,\eps},
\end{align}
and let $\mu_{\mathcal{G}}$ be a uniform probability distribution on the set $\mathcal{G}$, \ie, $\mu_\mathcal{G}(U_{\beta}) = \frac{1}{2^{2N}}$ for any $\beta$.

Then by definition \eqref{eqn::alg_error},
\begin{align}
\label{eqn::err_1}
\begin{aligned}
e_{\mathcal{F}, \Lambda}(A)^2 & \ge \ee_{\revise{U_{\beta}\sim\mu_\mathcal{G}}}\ee_{(\omega, \wt{\omega})\sim \prob\times \wt{\prob}}\Big[ \abs{X_T(U_{\beta}, \omega) - A(U_{\beta}, \omega, \wt{\omega})}^2 \Big]\\
&= \ee_{(\omega, \wt{\omega})\sim \prob\times \wt{\prob}}\bigg[ \ee_{\revise{U_{\beta}\sim\mu_\mathcal{G}}} \Big[ \abs{X_T(U_{\beta}, \omega) - A(U_{\beta}, \omega, \wt{\omega})}^2\Big]\bigg] \\
&\ge \prob(\mathcal{E})\ \ee_{(\omega, \wt{\omega}) \sim  \prob\rvert_{\mathcal{E}}\times \wt{\prob}} \bigg[ \ee_{\revise{U_{\beta}\sim\mu_\mathcal{G}}} \Big[ \abs{X_T(U_{\beta}, \omega) - A(U_{\beta}, \omega, \wt{\omega})}^2 \Big] \bigg],
\end{aligned}
\end{align}
\revise{where the probability measure $\prob\rvert_{\mathcal{E}}$ is the restriction of
$\mathbb{P}$ to the event $\mathcal{E}$; more specifically,
$\prob\rvert_{\mathcal{E}} (B) = \prob(B\cap \mathcal{E})/\prob(\mathcal{E})$
for any event $B\in \Sigma$.
Besides, $\ee_{{U_{\beta}\sim\mu_\mathcal{G}}}$ means expectation with respect to $U_{\beta}$, where $U_{\beta}$ is treated as a uniformly distributed random variable on the set $\mathcal{G}$; to avoid introducing too many notations, we slightly abuse the notation of $U_{\beta}$ to represent both a single potential function (or say an element in $\mathcal{G}$) and a (uniformly distributed) random variable on the set $\mathcal{G}$.
}

We know \revise{by \lemref{lem::choose_bounds}} that $\prob(\mathcal{E})$ is uniformly bounded \revise{from below} by a positive value $\mathsf{P}(C_x, C_v, u, T)$, and we claim (to be proved below in the \emph{Step (2)})
\begin{align}
\label{eqn::err_2}
\ee_{\revise{U_{\beta}\sim\mu_\mathcal{G}}} \Big[ \abs{X_T(U_{\beta}, \omega) - A(U_{\beta}, \omega,\wt{\omega})}^2\Big]  \gtrsim\ \frac{C_x^4 \overline{C}^2 \xi^2}{N^3}.\end{align}
Therefore,
\begin{align*}
e_{\mathcal{F}, \Lambda}(A) \gtrsim C_{\text{low}} N^{-3/2},
\end{align*}
where
\begin{align*}
\begin{aligned}
C_{\text{low}} &= \sqrt{\mathsf{P}(C_x,  C_v, u, T)} C_x^2 \overline{C} \xi \\
&\myeq{\eqref{eqn::parameters}} \sqrt{\mathsf{P}(C_x,  C_v, u, T)} C_x^2 \overline{C} \min\{ u-\ell, u_R - u\}.
\end{aligned}
\end{align*}
Since $C_x$, $C_v$, $u$, and $u_R$ are parameters to tune, we arrive at \eqref{eqn::c_low} by optimizing.

\medskip

{\noindent \emph{Step ({2}): Proof of the perturbation bound \eqref{eqn::err_2}.}}
From now on, we fix both $\wt{\omega}\in \wt{\probsp}$ and $\omega\in \mathcal{E}$. The main task is to estimate the fluctuation of the exact solution $X_T(U_{\beta}, \omega)$ for those $U_{\beta}$ with the same algorithmic output  $A(U_{\beta}, \omega, \wt{\omega})$ by \propref{prop::diff_x_lower}; this is done by a quantitative perturbation analysis.

The first task is to characterize the set of $U_{\beta}$ with the same algorithmic output $A(U_{\beta}, \omega, \wt{\omega})$. This is given by Lemma~\ref{lem::same_result} below. To state the result, let us define some notations.
The access points for $\nabla U_{\beta}$ are denoted by $Y_1^{\beta}, Y_2^{\beta}, \cdots, Y_N^{\beta}$, which only depend on the choice of $\beta$. Let us denote the union of sub-intervals that $Y_j^{\beta}$ belong to as $\mathcal{J}^{\beta}$:
\begin{align*}
    \revise{\mathcal{J}^{\beta} := \bigcup_{\substack{-N\le j \le N-1,\\ Y_k^{\beta} \in I_j \text{ for some } 1\le k\le N}} I_j.}
\end{align*}
If $\mu(\mathcal{J}^{\beta}) < C_x/2$ (\revise{\eg,} there exist two indices $j_1 < j_2$ such that $Y_{j_1}^{\beta}$ and $Y_{j_2}^{\beta}$ belong to the same sub-interval), then we add sub-intervals with the largest indices to complete $\mathcal{J}^{\beta}$.
\revise{
More specifically, for each fixed $\beta$, let us define 
\begin{align*}
\mathpzc{L} = \Big\{j\ \rvert\ -N\le j \le N-1,\ Y_{k}^{\beta} \notin I_j,\ \forall k\in \{1, 2, 3, \cdots, N\}\ \Big\}.
\end{align*}
When the cardinality of $\mathpzc{L}$ is larger than $N$ (i.e., $\abs{\mathpzc{L}} > N$), or equivalently $\mu(\mathcal{J}^{\beta}) < C_x/2$, we add $|\mathpzc{L}|-N$ sub-intervals in a procedure described by the following pseudocode:
\begin{center}
\begin{minipage}{.6\linewidth}
\floatname{algorithm}{Pseudocode: }
\begin{algorithm}[H]
\renewcommand\thealgorithm{}
\begin{algorithmic}
\WHILE {$\mu(\mathcal{J}^{\beta}) < C_x/2$}
\STATE $j\leftarrow $ largest value in $\mathpzc{L}$
\STATE $\mathpzc{L}\leftarrow \mathpzc{L} \backslash \{j\}$
\STATE $\mathcal{J}^{\beta} \leftarrow \mathcal{J}^{\beta} \cup I_j$
\ENDWHILE
\end{algorithmic}
\addtocounter{algorithm}{-1}
\caption{The completion procedure of $\mathcal{J}^{\beta}$}
\end{algorithm}
\end{minipage}
\end{center}
}

\smallskip
\begin{example}
\normalfont
If $N = 3$, then there are six sub-intervals $I_{-3}, I_{-2}, I_{-1}, I_{0}, I_{1}, I_{2}$. Let us consider the following examples:
\begin{itemize}
\item  if all $Y_1^{\beta}, Y_2^{\beta}, Y_3^{\beta} \in I_{0}$, then we set $\mathcal{J}^{\beta} := I_{0} \cup (I_{1}\cup I_{2})$;
\item if all $Y_1^{\beta}, Y_2^{\beta}, Y_3^{\beta} \in I_{-1}$, then we set $\mathcal{J}^{\beta} := I_{-1} \cup (I_{1}\cup I_{2})$;
\item if $Y_1^{\beta}, Y_2^{\beta} \in I_{0}$, and $Y_3^{\beta}\in I_1$, then we set $\mathcal{J}^{\beta} := I_{0} \cup I_{1} \cup (I_2)$;
\item If $Y_1^{\beta}\in I_0$, $Y_2^{\beta}\in I_1$, and $Y_3^{\beta} \in I_2$, then we set $\mathcal{J}^{\beta} := I_{0} \cup I_{1} \cup I_2$.
\end{itemize}
The interval within the parenthesis is the additional sub-intervals that we add to complete $\mathcal{J}^{\beta}$.
\end{example}
\smallskip

Note that such a procedure is always valid, since $Y_1^{\beta}, Y_2^{\beta}, \cdots, Y_N^{\beta}$ reside in at most $N$ sub-intervals; for each $\beta$, $\mathcal{J}^{\beta}$ is always uniquely defined.
From now on, when we use the notation $\mathcal{J}^{\beta}$, we refer to the \enquote{completed} version.
Likewise, time points for queries, denoted by $t_j^{\beta}$,  also depend on $\beta$ only.

\smallskip
\begin{lemma}
\label{lem::same_result}
For any fixed index $\beta$, if for some index $\beta'$, \revise{$\nabla U_{\beta'} = \nabla U_{\beta}$ on the domain $\mathcal{J}^{\beta}$}, then
\begin{enumerate}[(i)]
\item $Y_j^{\beta} = Y_j^{\beta'}$, $t_j^{\beta} = t_j^{\beta'}$ for all $1\le j\le N$;
\item $\mathcal{J}^{\beta} = \mathcal{J}^{\beta'}$ and moreover,
$
A(U_{\beta}, \omega, \wt{\omega}) = A(U_{\beta'}, \omega, \wt{\omega}).
$
\end{enumerate}
\end{lemma}
\smallskip

\begin{proof}
Part (ii) trivially follows from the form of algorithms in \eqref{eqn::alg_form} and part (i).
Then it suffices to prove part (i). This comes from induction: if  $Y_j^{\beta} = Y_j^{\beta'}$, $t_j^{\beta} = t_j^{\beta'}$ for all $j \le k$, then we know the information up to $k$ queries are the same, \ie{}, $\Phi_k^{\beta} = \Phi_k^{\beta'}$ where $\Phi_k$ is defined in \eqref{eqn::info_j}, and superscripts are used to indicate the dependence;
\revise{to verify $\Phi_k^{\beta} = \Phi_{k}^{\beta'}$, we need to check the following two relations for any $j\le k$:
\begin{align*}
    \bmpair_{t_j^{\beta}}(\omega) = \bmpair_{t_j^{\beta'}}(\omega),\ 
    \nabla U_{\beta'}(Y_j^{\beta'}) = \nabla U_{\beta}(Y_j^{\beta}).
\end{align*}
Since $t_j^{\beta} = t_j^{\beta'}$ by the assumption of induction, obviously, $\bmpair_{t_j^{\beta}}(\omega) = \bmpair_{t_j^{\beta'}}(\omega)$.
Next, since $Y_j^{\beta} = Y_j^{\beta'}$, we have $\nabla U_{\beta'}(Y_{j}^{\beta'}) = \nabla U_{\beta'}(Y_j^\beta)$. It is then enough to verify $\nabla U_{\beta'}(Y_j^{\beta}) = \nabla U_{\beta}(Y_j^{\beta})$.
\begin{itemize}
\item If $Y_j^{\beta}\notin \bigl[-C_x/2, C_x/2\bigr]$,  then $\nabla U_{\beta'}(Y_j^{\beta})\ \myeq{\eqref{eqn::U_family}}\ u Y_{j}^{\beta}\ \myeq{\eqref{eqn::U_family}}\  \nabla U_{\beta}(Y_j^{\beta})$.

\item If $Y_j^{\beta}\in \bigl[-C_x/2, C_x/2\bigr]$, then $Y_j^{\beta}\in \mathcal{J}^{\beta}$ by definition. 
By the assumption that $\nabla U_{\beta'} = \nabla U_{\beta}$ on the domain $\mathcal{J}^{\beta}$, we know $\nabla U_{\beta'}(Y_j^{\beta}) = \nabla U_{\beta}(Y_j^{\beta})$.
\end{itemize}
}
Then
\revise{$(Y_{k+1}^{\beta}, t_{k+1}^{\beta}) = \varphi_k(Y_1^{\beta}, \Phi_k^{\beta}) = \varphi_k(Y_1^{\beta'}, \Phi_k^{\beta'}) = (Y_{k+1}^{\beta'}, t_{k+1}^{\beta'})$. 
As for the base case, note that for an arbitrary $\beta$, $U_{\beta}\in \mathcal{F}_{u,\eps}$, and thus $Y_1^{\beta}\ \myeq{\eqref{eqn::Y1_t1}}\ x^{\star}(U_{\beta})\  \myeq{\eqref{eqn::Fueps}}\ 0$.}
\hfill
\end{proof}

\smallskip
\begin{definition}
For two arbitrary indices $\beta$ and $\beta'$, we define a binary relation $\beta\sim \beta'$ if $\beta'_j = \beta_j$ whenever $I_j\subseteq \mathcal{J}^{\beta}$.
By the above \lemref{lem::same_result}, we know $\mathcal{J}^{\beta} = \mathcal{J}^{\beta'}$ and $A(U_{\beta}, \omega, \wt{\omega}) = A(U_{\beta'}, \omega, \wt{\omega})$. It is easy to verify that such a relation is an equivalence relation.
\end{definition}
\smallskip

For any index $\beta$, apparently, there are exactly $2^{N}-1$ other indices belonging to the same equivalence class (since we could freely choose $\beta_j \in \{0,1\}$ whenever $I_j \not\subseteq \mathcal{J}^{\beta}$),
and there are exactly $2^N$ such equivalence classes. Let us enumerate these equivalence classes by $\mathcal{K}_1, \mathcal{K}_2, \cdots, \mathcal{K}_{2^N}$.

\smallskip
\begin{example}
\normalfont
When $N = 3$, for an index $\beta = (0, 0, 0, 0, 0, 0)$, we assume that
$\mathcal{J}^{\beta} = I_0 \cup I_1 \cup I_2$ as an example. Then we could
freely choose the first three indices, and the equivalence class containing
$\beta$ is exactly \revise{$\Big\{(\beta_{-3},\beta_{-2},\beta_{-1}, 0, 0, 0)\
    \Big\rvert\
    \beta_{-3},\beta_{-2},\beta_{-1}\in \{0,1\}\Big\}$.}
\end{example}
\smallskip

We now consider how much the actual solution $X_T(U_{\beta}, \omega)$ can fluctuate within the same class.
Consider an arbitrary equivalence class $\mathcal{K}$, and suppose $\beta\in
\mathcal{K}$. For any index $\beta'$, recall that $\mathcal{J}^{\beta'}$ are
the same for all $\beta'\in \mathcal{K}$ 
(\ie{}, \revise{$\mathcal{J}^{\beta'}$} only depends on the equivalence class $\mathcal{K}$ that we consider), 
and $\beta'_j$ are the same if $I_j \subseteq \mathcal{J}^{\beta}$.
This motivates us to define the reduced index below.

\smallskip
\begin{definition}[Reduced index]
    For an equivalence class $\mathcal{K}$, \revise{let us denote $\mathcal{J}
    \equiv \mathcal{J}^{\beta}$ for an arbitrary $\beta\in \mathcal{K}$; from the explanation
    above, the set $\mathcal{J}$ is well-defined.}
We introduce the \emph{reduced index} $\wt{\beta} := \big(\beta_{j}\big)_{I_j \not\subseteq \mathcal{J}} \in \{0,1\}^{N}$.
For each class $\mathcal{K}$, there is a one-to-one correspondence between $\beta\in \mathcal{K}$ and a reduced index $\wt{\beta}\in\{0,1\}^{N}$. Thus, we slightly abuse the notation and denote $U_{\wt{\beta}}\equiv U_{\beta}$, whenever the equivalence class $\mathcal{K}$ is clear from the context.
For the reduced index, we define a \emph{partial order $\succ$} as follows: if $\wt{\beta}'_j \ge \wt{\beta}_j$ for all $1\le j \le N$ (namely, $\nabla U_{\beta'}(x) \ge \nabla U_{\beta}(x)$ for all $x\in \Real$), then we denote $\wt{\beta}' \succ \wt{\beta}$ (or $\wt{\beta}\prec \wt{\beta}'$).
\end{definition}
 \smallskip

\revise{Recall from the Step (1) above that $\mathscr{I} = \big[\frac{C_x}{8N}, \frac{3C_x}{8N}\big]$ with length $\mu(\mathscr{I}) = \frac{C_x}{4N}$.}
By \propref{prop::diff_x_lower} with $\mathcal{I} = \bigcup_{\substack{j:\ \wt{\beta}'_j > \wt{\beta}_j}} (x_j + \mathscr{I})$, we immediately have the following result.

\smallskip
\begin{lemma}
If $\wt{\beta}' \succ \wt{\beta}$, then
\begin{align}
\label{eqn::err_3}
\begin{aligned}
\abs{X_T(U_{\beta'}, \omega) - X_T(U_{\beta}, \omega)} &\ge  \overline{C} \eps \frac{C_x}{4N} \#\big\{j:\ \wt{\beta}'_j > \wt{\beta}_j \big\} \\
&= \frac{C_x^2 \overline{C} \xi}{32 N^2} \#\big\{j:\ \wt{\beta}'_j > \wt{\beta}_j \big\}.
\end{aligned}
\end{align}
\end{lemma}
\smallskip

Let us introduce the following set, for any integer $0\le k\le N$,
\begin{align*}
\mathcal{M}_k := \Big\{\wt{\beta}\in\{0,1\}^{N}:\ \sum_{j=1}^{N} \wt{\beta}_j = k \Big\}.
\end{align*}

\smallskip
\begin{lemma}
\label{lem::counting_order}
Consider even integer $N \ge 2$.
For any integer $0\le k \le N/2$, there exists a bijective map $\Upsilon: \mathcal{M}_k \rightarrow \mathcal{M}_{N-k}$ such that for any
$\wt{\beta}\in \mathcal{M}_k$,
we have
$
\wt{\beta}\prec \Upsilon(\wt{\beta})
$.
In particular, when $k = N/2$, $\Upsilon(\wt{\beta}) = \wt{\beta}$.
\end{lemma}
\smallskip

\begin{proof}
This lemma follows immediately from the symmetric chain decomposition (SCD) for Boolean lattices; see, \eg{}, \cite{gregor_hypercube_nodate,zhu_boolean_nodate} for an introduction, as well as proofs.
A symmetric chain is a sequence $\gamma^{(n)}\prec \gamma^{(n+1)} \prec \cdots \prec \gamma^{(N-n)}$ where $\gamma^{(j)} \in \mathcal{M}_j$ for $n\le j \le N-n$.
SCD states that the set $\{0,1\}^N$ can be decomposed into disjoint symmetric chains. Therefore, for any $\wt{\beta}\in \mathcal{M}_k$, it must belong to a particular chain, say $\gamma^{(n)}\prec \gamma^{(n+1)} \prec \cdots \prec \gamma^{(N-n)}$.
By the definition of symmetric chains, we have $\gamma^{(k)} = \wt{\beta}$,
and then we can simply define $\Upsilon(\wt{\beta}) := \gamma^{(N-k)}\succ \gamma^{(k)} \equiv \wt{\beta}$.
Such a procedure is always valid, and since all symmetric chains are disjoint, $\Upsilon$ is a bijective map.
\hfill
\end{proof}
\smallskip

With these preparations, we can now continue to finish the proof of \eqref{eqn::err_2}, and thus the lower bound in \thmref{thm::underdamped_complexity_finite_nablaU} for the case $d=1$.
Within any equivalence class $\mathcal{K}$, we have
\begin{align*}
\sum_{\beta\in \mathcal{K}}  &\abs{X_T(U_{\beta}, \omega)  - A(U_{\beta}, \omega,\wt{\omega})}^2 = \revise{\sum_{\wt{\beta}\in \{0,1\}^N}  \abs{X_T(U_{\wt{\beta}}, \omega)  - A(U_{\wt{\beta}}, \omega,\wt{\omega})}^2} \\
\ge & \sum_{k=0}^{\frac{N}{2} - 1} \sum_{\wt{\beta}\in M_k} \abs{X_T(U_{\wt{\beta}}, \omega) - A(U_{\wt{\beta}}, \omega, \wt{\omega})}^2 + \abs{X_T(U_{\Upsilon(\wt{\beta})}, \omega) - A(U_{\Upsilon(\wt{\beta})}, \omega, \wt{\omega})}^2 \\
\gtrsim& \sum_{k=0}^{\frac{N}{2} - 1} \sum_{\wt{\beta}\in M_k} \abs{X_T(U_{\wt{\beta}}, \omega) - X_T(U_{\Upsilon(\wt{\beta})}, \omega)}^2\ 
\mygesim{\eqref{eqn::err_3}}\  \sum_{k=0}^{\frac{N}{2}-1} \binom{N}{k} \frac{C_x^4 \overline{C}^2 \xi^2}{N^4} (N - 2 k)^2\\
\gtrsim & \frac{C_x^4 \overline{C}^2 \xi^2}{N^4} \sum_{k=0}^{\frac{N}{2}} \binom{N}{k} \big(\frac{N}{2}-k\big)^2
\gtrsim  \frac{C_x^4 \overline{C}^2 \xi^2}{N^4} \sum_{k=0}^{N} \binom{N}{k} \big(\frac{N}{2} - k\big)^2
\gtrsim  \frac{C_x^4 \overline{C}^2 \xi^2}{N^3} 2^N.
\end{align*}
Finally, we have
\begin{multline*}
\ee_{\revise{U_{\beta}\sim\mu_\mathcal{G}}} \Big[ \abs{X_T(U_{\beta}, \omega) - A(U_{\beta}, \omega, \wt{\omega})}^2\Big]
=\ \frac{1}{2^{2N}} \sum_{\mathcal{K}_j}  \sum_{\beta\in \mathcal{K}_j} \abs{X_T(U_{\beta}, \omega) - A(U_{\beta}, \omega, \wt{\omega})}^2 \\
\gtrsim\ \frac{1}{2^{2N}} \times 2^N \times \frac{C_x^4 \overline{C}^2 \xi^2}{N^3} 2^N
\gtrsim\ \frac{C_x^4 \overline{C}^2 \xi^2}{N^3}.
\end{multline*}
Thus we complete the proof of the lower bound in \thmref{thm::underdamped_complexity_finite_nablaU} for the case $d=1$.

\medskip

{\noindent \bf Case (\rom{2}): General dimension $d$.}
For a general dimension $d$, we can choose \revise{$U(x) = U(x_1, x_2, \cdots, x_d) = \sum_{j=1}^d U^{(j)}(x_j)$}, where $U^{(j)}:\Real\rightarrow\Real$, \revise{$x = (x_1, x_2, \cdots, x_d)$, and $x_j\in \Real$}. If $U^{(j)}\in\mathcal{F}(1, \ell, L)$, then $U\in \mathcal{F}(d, \ell, L)$.  Note that if $U$ takes this particular form, then each component of the underdamped Langevin dynamics \eqref{eqn::underdamped} is evolving independently.
Then we immediately know
\begin{align*}
e_{\mathcal{F},\Lambda}(A)^2 &\ge \revise{\sup_{U:\  U(x) = \sum_{j=1}^{d} U^{(j)}(x_j)}} 
\ee_{\revise{(\omega, \wt{\omega})\sim \prob\times \wt{\prob}}}\big[ \abs{X_T(U, \omega) - A(U, \omega, \wt{\omega})}^2\big] \\
&= \revise{\sup_{U:\ U(x) = \sum_{j=1}^{d} U^{(j)}(x_j)}}  \sum_{j=1}^{d}\ \ee_{\revise{(\omega, \wt{\omega})\sim \prob\times \wt{\prob}}} \Big[ \abs{X_T^{(j)}(U, \omega) - A^{(j)}(U, \omega, \wt{\omega})}^2 \Big] \\
&\gtrsim\  d C_{\text{low}}^2 N^{-3},
\end{align*}
where $X_T^{(j)}(U, \omega)\in \Real$ is the $j^{\text{th}}$ component of $X_T(U,\omega)$, and 
$A^{(j)}(U, \omega, \wt{\omega})$ is the $j^{\text{th}}$ component of the algorithmic prediction $A(U, \omega, \wt{\omega})$.
Therefore, we obtain $e_{\mathcal{F},\Lambda}(A) \gtrsim C_{\text{low}}\sqrt{d} N^{-3/2}$.
The last inequality above is intuitively reasonable, since $A^{(j)}(U, \omega, \wt{\omega})$ could be regarded as an algorithm of the $j^{\text{th}}$ component, and having queries to the information from independent components, like \revise{$\nabla_{x_k} U^{(k)}(x_k)$} ($k\neq j$), would not improve the algorithmic prediction for the $j^{\text{th}}$ component.

More rigorously, one could directly generalize the proof of the {Case (\rom{1})}. Below is a sketch of the only few technical differences. First, similar to \eqref{eqn::err_1},
\begin{align*}
e_{\mathcal{F}, \Lambda}(A)^2 \ge \sum_{j=1}^{d} \revise{\ee_{(\omega,\wt{\omega})\sim \prob\times \wt{\prob}}}\  \ee_{\revise{U\sim\mu_{\mathcal{G}_1}\times \mu_{\mathcal{G}_2} \times \cdots \times \mu_{\mathcal{G}_d}}} \Big[\abs{X_T^{(j)}(U, \omega) - A^{(j)}( U, \omega, \wt{\omega})}^2\Big],
\end{align*}
where $\mu_{\mathcal{G}_j}$ is a uniform measure of $U^{(j)} \in \mathcal{G}_j$, similar to \eqref{eqn::Gset}, for $1\le j\le d$.
Without loss of generality, we only need to consider the component $j=1$, and we need to show that
\begin{align*}
\revise{\ee_{(\omega,\wt{\omega})\sim \prob\times \wt{\prob}}\ \ee_{U\sim\mu_{\mathcal{G}_1}\times \mu_{\mathcal{G}_2}\times \cdots \times \mu_{\mathcal{G}_d}} } \Big[\abs{X_T^{(1)}(U, \omega) - A^{(1)}( U, \omega, \wt{\omega})}^2\Big] \gtrsim C_{\text{low}}^2 N^{-3}.
\end{align*}
We shall fix $\omega \in \mathcal{E}$, where $\mathcal{E}$ is now defined in the same way as \eqref{eqn::event_E} by considering the components $X_t^{(1)}(U,\omega)$ and $V_t^{(1)}(U,\omega)$ only. We shall also fix $\wt{\omega}\in \wt{\probsp}$ and $U^{(j)}$ for $j \ge 2$. Then it suffices to prove that (cf. \eqref{eqn::err_2})
\begin{align*}
\ee_{\revise{U^{(1)}\sim\mu_{\mathcal{G}_1}}} \Big[\abs{X_T^{(1)}(U^{\revise{(1)}}, \omega) - A^{(1)}( U^{\revise{(1)}}, \omega, \wt{\omega})}^2\Big] \gtrsim\ \frac{C_x^4 \overline{C}^2 \xi^2}{N^3}.
\end{align*}
The proof of this inequality is essentially the same as the Step (2) in the {Case (\rom{1})}. The only minor difference is that $\mathcal{J}^{\beta}$ is now defined as the completed union of sub-intervals where the first components of $Y_1^{\beta}, Y_2^{\beta}, \cdots, Y_N^{\beta}$ reside in.

\subsection{Proof of results in \secref{sec::perturbation}.}
\label{subsec::proof_perturbation_lemmas}
\smallskip

\begin{proof}{\bf (Proof of \lemref{lem::perturbation}).}\,
\revise{Consider any fixed $U\in \mathcal{F}_{u, \eps}$ and any fixed $\omega\in \probsp$.}
Let us introduce $g(x):= \nabla U(x) - u x$, $\Delta_{X, t}:= X_t(U, \omega) - X_t(U_u, \omega)$,
and $\Delta_{V,t} := V_t(U, \omega) - V_t(U_u, \omega)$.
By assumption, $\norm{g}_{\infty} \le \eps$ and $\Delta_{X,0} = \Delta_{V,0} = 0$.
By \eqref{eqn::underdamped},
it is straightforward to derive that
\begin{align*}
\ud \begin{bmatrix} \Delta_{X,t} \\ \Delta_{V,t} \end{bmatrix}  = H \begin{bmatrix} \Delta_{X,t} \\ \Delta_{V,t} \end{bmatrix}\ud t + \begin{bmatrix} 0 \\
- \frac{1}{L} g(X_t(U, \omega)) \end{bmatrix} \ud t,\  \qquad H = \begin{bmatrix} 0 & 1 \\ - \frac{u}{L} & -2 \end{bmatrix}.
\end{align*}
Then we could rewrite the above equation in the integral form,
\begin{align*}
\begin{bmatrix} \Delta_{X,t} \\ \Delta_{V,t} \end{bmatrix}  = \int_{0}^{t} e^{H(t-s)} \begin{bsmallmatrix} 0 \\ -\frac{g(X_s(U, \omega))}{L} \end{bsmallmatrix}\ud s.
\end{align*}
Hence, by introducing $g_s \equiv g(X_s(U, \omega))$, and recalling $\lambda_{\pm}$ from \eqref{eqn::lambda}, we have
\begin{align}
\label{eqn::diff_xv_v4}
\begin{aligned}
\Delta_{X,t} &= -\frac{1}{2\sqrt{L(L-u)}} \int_{0}^{t} g_s \bigl(e^{-(t-s)\lambda_{-}}-e^{-(t-s) \lambda_{+}}\bigr)\ud s,\\
\Delta_{V,t} &= -\frac{1}{2\sqrt{L(L-u)}} \int_{0}^{t} g_s \bigl(\lambda_{+} e^{-(t-s) \lambda_{+}} - \lambda_{-} e^{-(t-s)\lambda_{-}}\bigr) \ud s.
\end{aligned}
\end{align}
Since $\abs{g_s}\le \eps$, it is straightforward to obtain that
\begin{align*}
\abs{\Delta_{X,t}} &\le \frac{\eps}{2\sqrt{L(L-u)}}\int_{0}^{t} \abs{e^{-(t-s)\lambda_{-}} - e^{-(t-s)\lambda_{+}}}\ud s \\
&= \frac{\eps}{2\sqrt{L(L-u)}} \bigl(\frac{1-e^{-\lambda_{-} t}}{\lambda_{-}} - \frac{1-e^{-\lambda_{+} t}}{\lambda_{+}}\bigr) \\
& \le \frac{\eps}{2 \lambda_{-} \sqrt{L(L-u)}} = \frac{\eps}{2 L (1-\sqrt{1-\frac{u}{L}}) \sqrt{1-\frac{u}{L}}}.
\end{align*}
Similarly, for $\abs{\Delta_{V,t}}$, we have
\begin{align*}
\abs{\Delta_{V,t}} &\le \frac{\eps}{2\sqrt{L(L-u)}} \int_{0}^{t} \abs{\lambda_{+} e^{-(t-s)\lambda_{+}} - \lambda_{-} e^{-(t-s)\lambda_{-}}}\ud s \\
&\le \frac{\eps}{2\sqrt{L(L-u)}} \int_{0}^{t} \lambda_{+} e^{-(t-s)\lambda_{+}} + \lambda_{-} e^{-(t-s)\lambda_{-}}\ud s \\
&= \frac{\eps}{2\sqrt{L(L-u)}}  (1- e^{-\lambda_{+} t} + 1 - e^{-\lambda_{-} t})
\le \frac{\eps}{L\sqrt{1- \frac{u}{L}}}.
\end{align*}\hfill
\end{proof}

\medskip

\begin{proof}{\bf (Proof of \lemref{lem::choose_bounds}).}\,
For any $C_x, C_v>0$, let us pick any $0<\eps<\bar{\epsilon}$.
Notice that
\[\sup_{0\le t\le T} X_t(U_{u}, \omega) \ge 2 C_x \ge C_x + \frac{\eps}{2 L (1-\sqrt{1-\frac{u}{L}}) \sqrt{1-\frac{u}{L}}}
\]
implies that $\sup_{0\le t\le T} X_t(U, \omega) \ge C_x$ by \eqref{eqn::diff_xv_upper}, and likewise for the other two cases.
We have
\begin{align*}
&\begin{aligned} \prob\Big(\omega: \sup_{0\le t\le T} X_t(U, \omega) \ge C_x,\ \inf_{0\le t\le T} & X_t(U, \omega) \le - C_x,\\
& \sup_{0\le t\le T} \abs{V_t(U, \omega)} \le C_v,\  \forall U\in \mathcal{F}_{u, \eps} \Big) \\
 \end{aligned}\\
 & \begin{aligned}
  \ge \prob\Big(\omega: \sup_{0\le t\le T} X_t(U_{u}, \omega) \ge 2C_x,\  \inf_{0\le t\le T} & X_t(U_{u}, \omega) \le - 2C_x, \\
  & \sup_{0\le t\le T} \abs{V_t(U_{u}, \omega)} \le C_v/2\Big).
  \end{aligned}
 \end{align*}
Finally, recall the expression of $\mathsf{P}(C_x, C_v, u, T)$ from \eqref{eqn::P_lower_prob}.
\hfill
\end{proof}
\medskip

\begin{proof}{\bf (Proof of \propref{prop::diff_x_lower}).}\,
We shall fixed $\omega\in\mathcal{E}$ throughout this proof.
Let $\Delta_{X,t} := X_t(U_1,\omega) - X_t(U_2, \omega)$ and $\Delta_{V,t} := V_t(U_1,\omega) - V_t(U_2, \omega)$.  Then by \eqref{eqn::underdamped}, we have
\begin{align}
\label{eqn::diff_xv_v2}
\ud \begin{bmatrix} \Delta_{X,t} \\ \Delta_{V,t} \end{bmatrix}  = H_t \begin{bmatrix} \Delta_{X,t} \\ \Delta_{V,t} \end{bmatrix}\ud t + \begin{bmatrix} 0 \\
- \frac{1}{L} g(X_t(U_1, \omega)) \end{bmatrix} \ud t, \qquad H_t = \begin{bmatrix} 0 & 1 \\ - \frac{u_t}{L} & -2 \end{bmatrix},
\end{align}
with initial conditions $\Delta_{X,0} = \Delta_{V,0} = 0$, where \revise{$u_{\cdot}: t\rightarrow u_t$} is a continuous function of time such that
\begin{align*}
\nabla U_2(X_t(U_1, \omega)) - \nabla U_2(X_t(U_2, \omega)) = u_t \Delta_{X,t}.
\end{align*}
\revise{To see why $u_t$ is well-defined}, notice that $U_2\in \mathcal{F}$ is a $C^2(\Real)$ function (recall $\mathcal{F}$ in \eqref{eqn::problem}).
Then by the first-order Taylor's expansion,
\begin{align*}
\nabla U_2(X_t(U_1, \omega)) - \nabla U_2(X_t(U_2, \omega))  = \nabla^2 U_2(\vartheta) (X_t(U_1,\omega) - X_t(U_2, \omega)),
\end{align*}
for some value $\vartheta$ between $X_t(U_1, \omega)$ and $X_t(U_2,\omega)$; we simply let $u_t = \nabla^2 U_2(\vartheta)$.
Moreover, we could easily observe that $\ell \le u_t \le u_R \le L$.
For simplicity, we shall again denote
\begin{align*}
g_t \equiv g(X_t(U_1, \omega)),
\end{align*}
here and below.

Intuitively, the ODE dynamics \eqref{eqn::diff_xv_v2} consists of two parts: the contraction part and the source part.
Suppose $u_t \equiv u$ is independent of time, then it is easy to observe that $e^{H t}$ is a contraction operator for large enough time $t$, and the ODE dynamics \eqref{eqn::diff_xv_v2} with $g\equiv 0$ will convergence to the origin exponentially fast, for any initial condition;
the source part $\begin{bsmallmatrix} 0 \\ -\frac{g_t}{L}\end{bsmallmatrix}$ will drag the velocity term (\ie, $\Delta_{V,t}$) towards the negative direction. Under the assumption that $\omega\in \mathcal{E}$, the term $g_t$ takes non-zero value at least for a period of $\mu(\mathcal{I})/C_v$.

When $u_t = u$ for all $t$, \eqref{eqn::diff_xv_v2} has an explicit solution shown below, similar to \eqref{eqn::exact_soln} and \eqref{eqn::diff_xv_v4} above.

\smallskip
\begin{lemma} Suppose $u_t = u$ for all $t\in[0,T]$ in \eqref{eqn::diff_xv_v2}, then $\Delta_{X,T} = \mathcal{S}(u, T)$, where
\begin{align*}
& \mathcal{S}(u, T) := \left\{\begin{aligned}
& \int_{0}^{T} \frac{e^{(t-T)(1+\sqrt{1-u/L})} - e^{(t-T)(1-\sqrt{1-u/L})}}{2\sqrt{L(L-u)}} g_t\ud t,\ & \text{ if } \ell \le u <  L; \\
& \int_{0}^{T} \frac{(t-T) e^{t-T}}{L}g_t\ud t,  & \text{ if } u = L.\\
\end{aligned}\right.
\end{align*}
\end{lemma}\smallskip

Observe that $\Delta_{X,T} \le 0$ for both cases, which inspires us to propose the following general result.

\smallskip
\begin{lemma}
\label{lem::trapping_region}
The region characterized by $\Delta_{V,t}\le - \Delta_{X,t}$ and $\Delta_{X,t} \le 0$ forms a trapping region for the dynamics \eqref{eqn::diff_xv_v2}.
Therefore,
the quantity $\Delta_{X,t} \equiv X_t(U_1,\omega) - X_t(U_2,\omega) \le 0$ for any $t\in [0,T]$.
\end{lemma}
\smallskip

\begin{proof}
To prove that the region formed by $\Delta_{V,t} \le - \Delta_{X,t}$ and $\Delta_{X,t} \le 0$ is a trapping region for the ODE dynamics \eqref{eqn::diff_xv_v2}, we consider the following three cases at the boundary:
\begin{itemize}

\item ($\Delta_{X,t} = 0$ and $\Delta_{V,t} = 0$). We know $\frac{\ud}{\ud t} \Delta_{X,t} = 0$ and $\frac{\ud}{\ud t} \Delta_{V, t} = -\frac{1}{L} g_t \le 0$. Thus, the solution of $\Delta_{X,t}$ and $\Delta_{V,t}$ will not escape the trapping region.

\item ($\Delta_{X,t} = 0$ and $\Delta_{V,t} < 0$, \ie, the negative half line of the velocity-axis). We know
\begin{align*}
(\frac{\ud}{\ud t} \Delta_{X,t}, \frac{\ud}{\ud t} \Delta_{V,t}) \cdot (-1, 0) = \big(\Delta_{V,t}, -2 \Delta_{V,t} - \frac{g_t}{L}\big) \cdot (-1, 0) = -\Delta_{V,t}  > 0.
\end{align*}
Therefore, the solution  of $\Delta_{X,t}$ and $\Delta_{V,t}$ will not escape the trapping region from the negative half-line of the velocity-axis.

\item ($\Delta_{V,t} = - \Delta_{X,t}$ for $\Delta_{X,t} < 0$).
\begin{align*}
(\frac{\ud}{\ud t} \Delta_{X,t}, \frac{\ud}{\ud t} \Delta_{V,t}) \cdot (-1, -1)
&= \big(\Delta_{V,t}, -\frac{u_t}{L} \Delta_{X,t} - 2 \Delta_{V,t} - \frac{g_t}{L}\big) \cdot (-1, -1) \\
=&\ (u_t/L - 1) \Delta_{X,t} + \frac{g_t}{L} \ge  0.
\end{align*}
By summarizing the above three cases, we conclude that the region formed by $\Delta_{V,t} \le - \Delta_{X,t}$ and $\Delta_{X,t}\le 0$ is indeed a trapping region. Since $\Delta_{X,0} = \Delta_{V,0} = 0$, we know that $\Delta_{X,t} \le 0$ for any time $t\ge 0$.
\end{itemize}
\hfill
\end{proof}

By the above lemma (\ie{}, $\Delta_{X,t}\le 0$ for any time $t\in [0,T]$), we know that
\begin{align}
\label{eqn::diff_xv_v3}
\begin{split}
\frac{\ud}{\ud t} \Delta_{X,t} &= \Delta_{V,t}, \\
\frac{\ud}{\ud t} \Delta_{V,t} &= -\frac{u_t}{L} \Delta_{X,t} - 2 \Delta_{V,t} - \frac{g_t}{L} \le -\frac{u_R}{L}\Delta_{X,t} - 2\Delta_{V,t}  - \frac{g_t}{L}.
\end{split}
\end{align}

\smallskip
\begin{lemma}
We claim that in general,
\begin{align}
\label{eqn::diff_x_upper_bound_by_u}
\Delta_{X,T} \le \mathcal{S}(u_R, T).
\end{align}
\end{lemma}
\smallskip

\begin{proof}
Let $(\Delta_{X,t}^{(1)}, \Delta_{V,t}^{(1)})$  \revise{be the solution of
\eqref{eqn::diff_xv_v2},} 
and let $(\Delta_{X,t}^{(2)}, \Delta_{V,t}^{(2)})$ be the solution of \eqref{eqn::diff_xv_v2} for $u_t \equiv u_R$.
Then let us introduce $\Gamma_{X,t} := \Delta_{X,t}^{(1)} - \Delta_{X,t}^{(2)}$
and $\Gamma_{V,t} := \Delta_{V,t}^{(1)} - \Delta_{V,t}^{(2)}$. 
\revise{By \eqref{eqn::diff_xv_v3}, we immediately have
\begin{align*}
    \frac{\ud}{\ud t}{\Gamma}_{X,t} &= \Gamma_{V,t},  \qquad
    \frac{\ud}{\ud t}{\Gamma}_{V,t} \le - \frac{u_R}{L}\Gamma_{X,t} - 2 \Gamma_{V,t}.
\end{align*}}

Since $\Gamma_{X,0} = \Gamma_{V,0} = 0$, by the same argument as in \lemref{lem::trapping_region}, we know that $\Gamma_{X,t}\le 0$ for any $t \ge 0$. 
Therefore, \revise{$\Delta_{X,T}^{(1)} \le \Delta_{X,T}^{(2)} \equiv
\mathcal{S}(u_R, T)$.} 
\hfill
\end{proof}

Let us consider the case $u_R < L$.
Then we have
 \begin{align*}
&\abs{\Delta_{X,T}} = -\Delta_{X,T} \\
& \myge{\eqref{eqn::diff_x_upper_bound_by_u}}\ \int_{0}^{T} \frac{e^{(t-T)(1-\sqrt{1-u_R/L})} - e^{(t-T)(1+\sqrt{1-u_R/L})}}{2\sqrt{L(L-u_R)}} g_t\ud t\\
&\begin{aligned}\myge{\eqref{eqn::assump_g}}\ \frac{\eps}{4 L\sqrt{1-\frac{u_R}{L}}} \int_{0}^{T} &e^{(t-T)(1-\sqrt{1-u_R/L})} \times\\
&\Big( 1 - e^{2(t-T)\sqrt{1-u_R/L}} \Big)\unit_{\mathcal{I}}(X_t(U_1,\omega))\ \ud t.\end{aligned}
\end{align*}
Without loss of generality, we assume $\tau_1(\omega) < \tau_2(\omega)$, where $\tau_1$ is the first hitting time of $X_t(U_1, \omega)$ to $-C_x$ and $\tau_2$ is the first hitting time to $C_x$. Since we assume $\omega\in \mathcal{E}$, both $\tau_1$ and $\tau_2$ are well-defined and $0\le \tau_1,\tau_2\le T$.
Then
 \begin{align*}
 &\abs{\Delta_{X,T}} \\
 &\ge \frac{\eps}{4 L\sqrt{1-\frac{u_R}{L}}} \int_{\tau_1}^{\tau_2} e^{(t-T)(1-\sqrt{1-\frac{u_R}{L}})} \Big( 1 - e^{2(t-T)\sqrt{1-\frac{u_R}{L}}} \Big)\unit_{\mathcal{I}}(X_t(U_1,\omega))\ \ud t \\
&= \frac{\eps}{4 L\sqrt{1-\frac{u_R}{L}}} \int_{\tau_1+\frac{C_x}{2C_v}}^{\tau_2 - \frac{C_x}{2C_v}} e^{(t-T)(1-\sqrt{1-\frac{u_R}{L}})} \Big( 1 - e^{2(t-T)\sqrt{1-\frac{u_R}{L}}} \Big)\unit_{\mathcal{I}}(X_t(U_1,\omega))\ \ud t \\
&\ge \frac{\eps \big(1 - e^{2(\tau_2 - \frac{C_x}{2C_v}-T)\sqrt{1-\frac{u_R}{L}}}\big)}{4 L\sqrt{1-\frac{u_R}{L}}} \int_{\tau_1+\frac{C_x}{2C_v}}^{\tau_2 - \frac{C_x}{2C_v}} e^{(t-T)(1-\sqrt{1-\frac{u_R}{L}})} \unit_{\mathcal{I}}(X_t(U_1,\omega))\ \ud t \\
&\ge \frac{\eps \big(1 - e^{- \frac{C_x}{C_v}\sqrt{1-\frac{u_R}{L}}}\big)}{4 L\sqrt{1-\frac{u_R}{L}}} \int_{\tau_1+\frac{C_x}{2C_v}}^{\tau_2 - \frac{C_x}{2C_v}} e^{(t-T)(1-\sqrt{1-\frac{u_R}{L}})} \unit_{\mathcal{I}}(X_t(U_1,\omega))\ \ud t,
 \end{align*}
where we use the following observation in the second equality:
since the velocity $V_t(U_1, \omega)$ is bounded by $C_v$, starting from the time $\tau_1$ (note that $X_{\tau_1}(U_1, \omega) = -C_x$), it takes at least $C_x/(2C_v)$ amount of time to reach $-C_x/2$. Thus, we have $\unit_{\mathcal{I}}(X_t(U_1, \omega)) = 0$ for $t\in \big[\tau_1, \tau_1 + \frac{C_x}{2C_v}\big]$ by the assumption \eqref{assump::support_I} that $\mathcal{I}\subseteq [-C_x/2, C_x/2]$ (similarly for the time period $\big[\tau_2 - \frac{C_x}{2 C_v}, \tau_2\big]$).

By the fact that $\mu\big(t: X_t(U_1,\omega) \in \mathcal{I}\big) \ge \frac{\mu(\mathcal{I})}{C_v}$ (namely, $\unit_{\mathcal{I}}(X_t(U_1,\omega)) = 1$ for at least $\mu(\mathcal{I})/C_v$ period of  time), we have
\begin{align*}
 &\abs{\Delta_{X,T}}  \\
 &\ge \frac{\eps \big(1 - e^{- \frac{C_x}{C_v}\sqrt{1-\frac{u_R}{L}}}\big)}{4 L\sqrt{1-\frac{u_R}{L}}} \int_{\substack{[\tau_1+\frac{C_x}{2C_v}, \tau_2 - \frac{C_x}{2C_v}] \cap \big\{t:\ X_t(U_1,\omega)\in \mathcal{I}\big\}}} e^{(t-T)(1-\sqrt{1-\frac{u_R}{L}})}\ \ud t \\
 &\ge \frac{\eps \big(1 - e^{- \frac{C_x}{C_v}\sqrt{1-\frac{u_R}{L}}}\big)}{4 L\sqrt{1-\frac{u_R}{L}}} \int_{\tau_1 + \frac{C_x}{2 C_v}}^{\tau_1 + \frac{C_x}{2 C_v} + \frac{\mu(\mathcal{I})}{C_v}} e^{(t-T)(1-\sqrt{1-\frac{u_R}{L}})}\ \ud t\\
 &= \frac{\eps \big(1 - e^{- \frac{C_x}{C_v}\sqrt{1-\frac{u_R}{L}}}\big)}{4 L\sqrt{1-\frac{u_R}{L}}} e^{(\tau_1 + \frac{C_x}{2C_v} - T) (1-\sqrt{1-\frac{u_R}{L}})} \frac{e^{\frac{\mu(\mathcal{I})}{C_v} (1-\sqrt{1-\frac{u_R}{L}})} - 1}{1-\sqrt{1-\frac{u_R}{L}}} \\
 &\ge \frac{\eps \big(1 - e^{- \frac{C_x}{C_v}\sqrt{1-\frac{u_R}{L}}}\big)}{4 L\sqrt{1-\frac{u_R}{L}}} e^{(\tau_1 + \frac{C_x}{2C_v} - T) (1-\sqrt{1-\frac{u_R}{L}})} \frac{\mu(\mathcal{I})}{C_v} \\
 &\ge \frac{\eps \mu(\mathcal{I}) \big(1 - e^{- \frac{C_x}{C_v}\sqrt{1-\frac{u_R}{L}}}\big)}{4 LC_v\sqrt{1-\frac{u_R}{L}}} e^{(\frac{3C_x}{2C_v} - T) (1-\sqrt{1-\frac{u_R}{L}})},
\end{align*}
where to get the final inequality, we use the fact that $\tau_1 \ge
\frac{C_x}{C_v}$ by the same {velocity-type argument}. Then we finish the proof
for the case $u_R < L$. 
As for $u_R = L$, \revise{one could either follow the similar
approach like above, or simply pass the limit $u_R\rightarrow L$ to obtain the
expression of $\overline{C}$}.
\hfill
\end{proof}

\subsection{Proof of the upper bound estimate for \thmref{thm::underdamped_complexity_finite_nablaU}.}
\label{subsec::proof_upper}
The upper bound estimate is based on Ref. \cite{shen_randomized_2019}. In the following, we shall provide a sketch only.
We consider non-adaptive mesh grid, \ie{}, $\timermm_k = k h$ with $h = T/\stepn$ as in \cite{shen_randomized_2019}.

Suppose $\hat{X}_{k}$, $\hat{X}_{k+1/2}$ and $\hat{V}_{k}$ are given by the randomized midpoint method as in \eqref{eqn::randomized_alg}, 
and suppose \revise{$\wt{X}_{k+1}$ and $\wt{V}_{k+1}$ are solutions of
\eqref{eqn::underdamped} at time $\timermm_{k +1}$,
conditioned on $X_{s_k} = \hat{X}_k$ and $V_{s_k} = \hat{V}_k$ at time $\timermm_k$}.
By \cite[Appendix E]{shen_randomized_2019}, we have
\begin{align*}
& \ee\big[\norm{\hat{X}_\stepn - X_{\reviseT}}^2 + \norm{\hat{X}_\stepn + \hat{V}_{\stepn} - X_{\reviseT} - V_{\reviseT}}^2\big] \\
\le &\ e^{-\stepn \frac{h}{2\varkappa}} \ee\big[\norm{\hat{X}_0 - X_{0}}^2 + \norm{\hat{X}_0 + \hat{V}_{0} - X_{0} - V_{0}}^2\big] \\
& \revise{+ \frac{2\varkappa}{h} \sum_{k=0}^{\stepn-1} \Big(3 \ee \normbig{\ee_{{\alphaRD_{k}}}[\hat{X}_{k+1} - \wt{X}_{k+1}]}^2 + 2 \ee\normbig{\ee_{{\alphaRD_{k}}}[\hat{V}_{k+1} - \wt{V}_{k+1}]}^2 \Big)}\\
&\revise{+ \sum_{k=0}^{\stepn-1} \Big(3 \ee \norm{\hat{X}_{k+1} - \wt{X}_{k+1}}^2 + 2 \ee \norm{\hat{V}_{k+1} - \wt{V}_{k+1}}^2\Big)} \\
=&\ \revise{\frac{2\varkappa}{h} \sum_{k=0}^{\stepn-1} \Big(3 \ee \normbig{\ee_{{\alphaRD_{k}}}[\hat{X}_{k+1} - \wt{X}_{k+1}]}^2 + 2 \ee\normbig{\ee_{{\alphaRD_{k}}}[\hat{V}_{k+1} - \wt{V}_{k+1}]}^2 \Big)}\\
&\revise{+ \sum_{k=0}^{\stepn-1} \Big(3 \ee \norm{\hat{X}_{k+1} - \wt{X}_{k+1}}^2 + 2 \ee \norm{\hat{V}_{k+1} - \wt{V}_{k+1}}^2\Big)},
\end{align*}
where in the last step, we use the fact that $\hat{X}_0 = X_{0}$ and $\hat{V}_0 = V_{0}$. \revise{In the above, $\ee_{\eta_{k}}$ means expectation with respect to $\eta_{k}$ for $0\le k \le \stepn-1$.}

By \cite[Lemma 2]{shen_randomized_2019}, we have
\begin{align*}
& \ee\big[\norm{\hat{X}_\stepn - X_{\reviseT}}^2 + \norm{\hat{X}_\stepn + \hat{V}_{\stepn} - X_{\reviseT} - V_{\reviseT}}^2\big]\\
\lesssim&\  \frac{2\varkappa}{h} \bigl( h^8 \sum_{k=0}^{\stepn-1} \ee \norm{\hat{V}_k}^2 + \frac{h^{10}}{L^2} \sum_{k=0}^{\stepn-1} \ee \norm{\nabla U(\hat{X}_k)}^2 + \frac{\stepn d h^9}{L} \big) \\
&+ h^4 \sum_{k=0}^{\stepn-1} \ee \norm{\hat{V}_k}^2 + \frac{h^4}{L^2}  \sum_{k=0}^{\stepn-1} \ee \norm{\nabla U(\hat{X}_k)}^2 + \frac{\stepn d h^5}{L} \\
\lesssim &  (h^4 + h^7 \varkappa) \sum_{k=0}^{\stepn-1}\ee \norm{\hat{V}_k}^2 + (\frac{h^9}{\ell L} + \frac{h^4}{L^2})\sum_{k=0}^{\stepn-1} \ee \norm{\nabla U(\hat{X}_k)}^2 + \bigl(\frac{\stepn d h^8}{\revise{\ell}} + \frac{\stepn d h^5}{L}\bigr).
\end{align*}
Next, we use \cite[Lemma 12]{shen_randomized_2019}, and obtain that
\begin{align*}
&\  \ee\big[\norm{\hat{X}_\stepn - X_{\reviseT}}^2 + \norm{\hat{X}_\stepn + \hat{V}_{\stepn} - X_{\reviseT} - V_{\reviseT}}^2\big] \\
\lesssim& \ (h^4 + h^7 \varkappa) \Big(\frac{\stepn d}{L} + \frac{1}{L} \Abs{\ee\bigl[\inner{\nabla U(\hat{X}_\stepn)}{\hat{V}_\stepn}\bigr]} \Big) \\
&+ (\frac{h^9}{\ell L}  + \frac{h^4}{L^2}) \Big(\stepn L d + \frac{L}{h} \Abs{\ee\bigl[\inner{\nabla U(\hat{X}_\stepn)}{\hat{V}_\stepn}\bigr]} \Big) + \bigl(\frac{\stepn d h^8}{\revise{\ell}} + \frac{\stepn d h^5}{L}\bigr) \\
\lesssim &\ \frac{h^3}{L} \Abs{\ee\bigl[ \inner{\nabla U(\hat{X}_\stepn)}{\hat{V}_\stepn}\bigr]} + \frac{h^4 \stepn d}{L},
\end{align*}
where in the last step, we use the fact that we are working on the $L^2$ strong error estimate and $h$ is the small parameter herein.
Then we need to estimate $\Abs{\ee\bigl[\inner{\nabla U(\hat{X}_\stepn)}{\hat{V}_\stepn}\bigr]}$. Similar to \cite[Appendix E]{shen_randomized_2019},
\begin{align*}
& \Abs{\ee\bigl[\inner{\nabla U(\hat{X}_\stepn)}{\hat{V}_\stepn}\bigr]} \\
\lesssim&\ \revise{\frac{1}{L} \ee\bigl[\norm{\nabla U(\hat{X}_{\stepn})}^2\bigr] + L \ee\bigl[ \norm{\hat{V}_{\stepn}}^2\bigr]}\\
\lesssim&\ L \ee\bigl[\norm{\hat{V}_\stepn - V_{\reviseT}}^2 + \norm{\hat{X}_\stepn - X_{\reviseT}}^2\bigr]
+ L \ee[\norm{V_{\reviseT}}^2] + \frac{1}{L} \ee\bigl[\norm{\nabla U(X_{\reviseT})}^2\bigr] \\
&\begin{aligned}\lesssim &\ L \ee\big[\norm{\hat{X}_\stepn - X_{\reviseT}}^2 + \norm{\hat{X}_\stepn + \hat{V}_{\stepn} - X_{\reviseT} - V_{\reviseT}}^2\big] \\
&+ L \ee\bigl[\norm{V_{\reviseT}}^2\bigr] + \frac{1}{L} \ee\bigl[\norm{\nabla U(X_{\reviseT})}^2\bigr]. \end{aligned}
\end{align*}
By combining the last two equations,
\begin{align*}
\ee\bigl[\norm{\hat{X}_\stepn - X_{\reviseT}}^2\bigr] \le &\  \ee\big[\norm{\hat{X}_\stepn - X_{\reviseT}}^2 + \norm{\hat{X}_\stepn + \hat{V}_{\stepn} - X_{\reviseT} - V_{\reviseT}}^2\big] \\
\lesssim &\ \frac{h^3}{L} \Big(L \ee\bigl[\norm{V_{\reviseT}}^2\bigr] + \frac{1}{L} \ee\bigl[\norm{\nabla U(X_{\reviseT})}^2\bigr]\Big) + \frac{h^4 \stepn d}{L} \\
\lesssim &\ h^3 \Big( \ee \bigl[\norm{V_{\reviseT}}^2\bigr] + \ee\bigl[\norm{X_{\reviseT} - x^{\star}}^2\bigr]\Big) + \frac{h^4 \stepn d}{L}.
\end{align*}

Suppose $(Y_t, Z_t)$ is another solution of the underdamped Langevin dynamics \eqref{eqn::underdamped} with the initial distribution as  $\rho_{\infty}$. Then similar to \cite[Appendix E]{shen_randomized_2019},
\begin{align*}
 &\ \ee \bigl[\norm{V_{\reviseT}}^2\bigr] + \ee\bigl[\norm{X_{\reviseT} - x^{\star}}^2\bigr] \\
 \lesssim &\ \ee\bigl[\norm{V_{\reviseT} - Z_{\reviseT}}^2 + \norm{X_{\reviseT} - Y_{\reviseT}}^2\bigr] + \ee\bigl[\norm{Z_{\reviseT}}^2\bigr] + \ee\bigl[\norm{Y_{\reviseT} - x^{\star}}^2\bigr] \\
 \lesssim &\ \ee\bigl[\norm{V_{\reviseT} - Z_{\reviseT} + X_{\reviseT} - Y_{\reviseT}}^2 + \norm{X_{\reviseT} - Y_{\reviseT}}^2\bigr] + \frac{d}{L} + \frac{d}{\ell} \\
 \lesssim &\ e^{-\frac{T}{\varkappa}} \ee\bigl[\norm{V_{0} - Z_0 + X_0 - Y_0}^2 + \norm{X_0 - Y_0}^2\bigr]  + \frac{d}{\ell}
 \lesssim\ (e^{-\frac{T}{\varkappa}} + 1)\frac{d}{\ell} \lesssim  \frac{d}{\ell}.
\end{align*}

Finally, we have
\begin{align*}
\ee\bigl[\norm{\hat{X}_\stepn - X_{\reviseT}}^2\bigr] & \lesssim \frac{d h^3}{\ell} + \frac{h^4 \stepn d}{L}
\lesssim  \frac{d}{\stepn^3} \bigl(\frac{T^3}{\ell} + \frac{T^4}{L}\bigr).
\end{align*}
Thus, if the algorithm $A$ is the randomized midpoint method \eqref{eqn::randomized_alg},
\begin{align*}
\revise{e_{\mathcal{F}, \Lambda}(A) \lesssim C_{\text{up}} \sqrt{d} \stepn^{-3/2} \lesssim C_{\text{up}} \sqrt{d} N^{-3/2},}
\end{align*}
where $C_{\text{up}} = \sqrt{\frac{T^3}{\ell} + \frac{T^4}{L}}$, \revise{and recall that the number of queries is $N = 2\stepn$}.

\section*{Acknowledgments}
This work is supported in part by the National Science Foundation via grants DMS-1454939 and CCF-1934964 (Duke TRIPODS).
We would like to thank Sam Hopkins (University of Minnesota) for introducing the symmetric chain decomposition \cite{351643}.

\bibliographystyle{cms_prelim}
\bibliography{reference.bib}

\newpage
\appendix
\section{RMM in the framework of randomized algorithms}
\label{app::rmm-randomized}

Rewriting the RMM \eqref{eqn::randomized_alg} in the framework of randomized algorithms introduced in \secref{subsec::ibc} is not technically challenging but many details require much attention, especially notations. Therefore, we would like to present how the RMM fits into the framework of general randomized algorithms in details for readers' convenience, and also for the purpose of a rigorous proof and completeness.

Let us consider the RMM for $\stepn$ steps (time step $h = T/\stepn$), which means $N = 2 \stepn$ evaluations of $\nabla U$ and $\bmpair_t$. 
The probability space $\wt{\probsp} = [0,1]^{\stepn}$ and $\vec{\eta} = (\eta_0, \eta_2, \cdots, \eta_{\stepn-1}) \in \wt{\probsp}$ determines random time steps in \eqref{eqn::randomized_alg}; the $\sigma$-algebra $\wt{\Sigma}$ can be simply chosen as the Borel $\sigma$-algebra of $\wt{\probsp}$; the probability measure $\wt{\prob}$ corresponds to a uniform distribution on $[0,1]^{\stepn}$. Next, we need to rewrite the RMM with fixed $\vec{\eta} \in \wt{\probsp}$ as a deterministic algorithm; see \eqref{eqn::det_alg_maps} and \eqref{eqn::alg_form}. With fixed $\vec{\eta}$, the time step prediction is already complete. Therefore, we just need to explain those terms $\Phi_j$ and $Y_j$ in \eqref{eqn::det_alg_maps}, as well as the forms of $\varphi_j$ and $\phi$.

\smallskip
{\noindent {\bf Step 1: Corresponding $Y_j$, $\Phi_j$ in \eqref{eqn::det_alg_maps}.}} 

Let us denote $\timermm_{k+1/2} = \timermm_k + \eta_k h \equiv (k+\eta_k) h$ for integer $0\le k \le \stepn-1$, and introduce artificial terms $\hat{V}_{k+1/2} := \hat{V}_k$ for convenience of notations below.
Let us introduce $t_j = \timermm_{j/2}$ and $Y_j = \hat{X}_{(j-1)/2}$ for integer $1\le j \le 2\stepn$. Recall the notation $\Phi_j$ from  \eqref{eqn::info_j} that 
\begin{align*}
\Phi_j := \Big(\Upsilon_{U,\omega}(Y_1, t_1), \Upsilon_{U,\omega}(Y_2, t_2), \cdots,  \Upsilon_{U,\omega}(Y_j, t_j)\Big).
\end{align*}

\begin{center}
\begin{table}[h!]
\begin{center}
\begin{tabular}{cc}
\toprule
General framework in \secref{subsec::ibc} & RMM \eqref{eqn::randomized_alg}\\
\toprule
$N=2\stepn$ evaluations & $\stepn$ steps\\
$t_j$ & $\timermm_{j/2}$\\
$Y_j$ & $\hat{X}_{(j-1)/2}$\\
\toprule
\end{tabular}
\end{center}
\caption{This table summarizes the relationship between two sets of notations for the general framework and the RMM.}
\end{table}
\end{center}

\smallskip
{\noindent {\bf Step 2: Choice of $\varphi_j$ and $\phi$.}} In this part, we will show that $Y_{j+1}\equiv \hat{X}_{j/2}$ is a function of $Y_1 \equiv x^{\star}(U)$ and $\Phi_j$ (\ie{}, to determine $\varphi_j$), and we will also show that $\hat{X}_{\stepn} = \phi(Y_1, \Phi_{2\stepn})$ for some deterministic function $\phi$.

\smallskip

\begin{prop}
We can find deterministic mappings $\phimap{1}{j}$ and $\phimap{2}{j}$, for integer $1\le j \le 2\stepn$, such that 
\begin{align*}
\hat{X}_{j/2} &= \phimap{1}{j}(Y_1, \Phi_{j});\\
\hat{V}_{j/2} &= \phimap{2}{j}(Y_1, \Phi_j).
\end{align*}
\end{prop}
\smallskip

As a consequence, the choices of $\varphi_j$ and $\phi$ are straightforward:
\begin{align*}
\varphi_j(Y_1, \Phi_j)&:= (\phimap{1}{j}(Y_1, \Phi_j), t_{j+1}),\ \text{ for } 1\le j \le 2\stepn -1 \equiv N-1;\\
\phi(Y_1, \Phi_{2\stepn}) & := \phimap{1}{2\stepn} (Y_1, \Phi_{2\stepn}) \equiv \hat{X}_{\stepn}.
\end{align*}

\begin{proof}
We prove it by mathematical induction. Suppose the result is true for any integer $1\le j\le 2k$, where $k$ is some integer ($1 \le k \le \stepn-1$), then we show the forms of $\phimap{1}{j}$ and $\phimap{2}{j}$ for $j = 2k+1$ and $2k +2$. That is,
both $\hat{X}_{k+1/2}$ and $\hat{V}_{k+1/2}$ are functions of $Y_1\equiv x^{\star}(U)$ and $\Phi_{2k+1}$; both $\hat{X}_{k+1}$ and $\hat{V}_{k+1}$ are functions of $Y_1$ and $\Phi_{2k+2}$.

\begin{itemize}

\item {\bf (Case $j = 2k+1$).} By definition, $\hat{V}_{k+1/2} \equiv \hat{V}_k = \phimap{2}{2k}(Y_1, \Phi_{2k})$. As $\Phi_{2k+1}$ contains more information than $\Phi_{2k}$, obviously, $\hat{V}_{k+1/2}$ is a deterministic function of $Y_1$ and $\Phi_{2k+1}$. Next, we consider $\hat{X}_{k+1/2}$. 
By \eqref{eqn::randomized_alg}, \eqref{eqn::weighted_bm} and by recalling $t_j \equiv \timermm_{j/2}$, 
\begin{align*}
&\begin{aligned}
\hat{X}_{k+1/2} &= \hat{X}_k +  \frac{1-e^{-2 \alphaRD_k h}}{2} \hat{V}_k + \frac{1}{\sqrt{L}} \Big(W_{t_{2k+1}} - W_{t_{2k}} - e^{-2 t_{2k+1}} \bigl(\wt{W}^{(2)}_{t_{2k+1}} - \wt{W}^{(2)}_{t_{2k}}\bigr)\Big) \\
&\qquad -\frac{1}{2L} \big(\int_{0}^{\alphaRD_k h} 1-e^{2(s-\alphaRD_k h)}\ud s\bigr) \nabla U(\hat{X}_k).
\end{aligned}
\end{align*}
Note that $\hat{X}_k = \phimap{1}{2k}(Y_1, \Phi_{2k})$, $\hat{V}_k = \phimap{2}{2k}(Y_1, \Phi_{2k})$ by assumption, and also note that 
\begin{align*}
\Phi_{2k+1} &= \Big(\Phi_{2k}, \Upsilon_{U,\omega}(Y_{2k+1}, t_{2k+1}) \Big),\\
\Upsilon_{U,\omega}(Y_{2k+1}, t_{2k+1}) &= \bigl(\nabla U(Y_{2k+1}), \bmpair_{t_{2k+1}}\bigr) \equiv \bigl(\nabla U(\hat{X}_{k}), \bmpair_{t_{2k+1}}\bigr),
\end{align*}
we immediately know that 
\begin{align*}
\hat{X}_{k+1/2} = \phimap{1}{2k}(Y_1, \Phi_{2k}) +  \frac{1-e^{-2 \alphaRD_k h}}{2} \phimap{2}{2k}(Y_1, \Phi_{2k})  + f_{2k+1}(\Phi_{2k+1}),
\end{align*}
for some linear function $f_{2k+1}$ on $\Phi_{2k+1}$; the expression of $f_{2k+1}$ should be clear from above. 
Then, we can conclude that 
$\hat{X}_{k+1/2} = \phimap{1}{2k+1} (Y_1, \Phi_{2k+1})$ for some mapping $\phimap{1}{2k+1}$, which is defined iteratively based on $\phimap{1}{2k}$, $\phimap{2}{2k}$, and $f_{2k+1}$.\smallskip

\item {\bf (Case $j = 2k + 2$).} Similarly, both $\hat{X}_{k+1}$ and $\hat{V}_{k+1}$ are functions of $\hat{X}_k$, $\hat{V}_k$, $\nabla U(\hat{X}_{k+1/2})$, $\bmpair_{\timermm_k}$ and $\bmpair_{\timermm_{k+1}}$, and the latter three terms are included in the information $\Phi_{2k+2}$. By the same argument, we know that both $\hat{X}_{k+1}$ and $\hat{V}_{k+1}$ are functions of the initial condition $Y_1\equiv x^{\star}(U)$ and $\Phi_{2k+2}$.

\item The base cases for $j = 1, 2$ can be similarly verified.
\end{itemize}
\hfill
\end{proof}

\section{Additional proofs}
\label{app::supplementary}
\smallskip

\begin{proof}{\bf (Proof of \lemref{lem::non-trivial-P}).}\,
In the limit $C_v\rightarrow\infty$,
\begin{align*}
&\hspace{1em} \lim_{C_v\rightarrow\infty} \mathsf{P}(C_x, C_v, u, T)\\
&\begin{aligned}= &\lim_{C_v\rightarrow\infty}\prob\Big(\omega:\ \sup_{0\le t\le T} X_t(U_u,\omega) \ge 2 C_x,\ \inf_{0\le t\le T} X_t(U_u, \omega)\le -2 C_x,\\ &\hspace{12em}\sup_{0\le t\le T} \abs{V_t(U_u,\omega)} \le C_v/2\Big)
\end{aligned}\\
&=\ \prob\Big(\omega:\ \sup_{0\le t\le T} X_t(U_u,\omega) \ge 2 C_x,\ \inf_{0\le t\le T} X_t(U_u, \omega)\le -2 C_x\Big).
\end{align*}
Note that $\mathsf{P}(C_x, C_v, u, T)$ is a monotonically increasing function with respect to $C_v$. Therefore, for sufficiently large $C_v$, we must have
\begin{align*}
    \mathsf{P}(C_x, C_v, u, T) \ge \frac{1}{2} \prob\Big(\omega:\ \sup_{0\le t\le T} X_t(U_u,\omega) \ge 2 C_x,\ \inf_{0\le t\le T} X_t(U_u, \omega)\le -2 C_x\Big).
\end{align*}
For finite $C_x$, due to the fluctuation of Brownian motion, the probability that the trajectory $\{X_t(U_u,\omega)\}_{0\le t\le T}$ crosses both levels $2C_x$ and $-2C_x$ is expected to be strictly positive. More rigorously, 
\begin{align*}
&\ \prob\Big(\omega:\ \sup_{0\le t\le T} X_t(U_u,\omega) \ge 2 C_x,\ \inf_{0\le t\le T} X_t(U_u, \omega)\le -2 C_x\Big) \\
\ge &\  \prob\Big(\omega:\ X_T(U_u,\omega) \ge 2 C_x,\ X_{T/2}(U_u, \omega) \le -2 C_x\Big).  
\end{align*}

By \eqref{eqn::exact_soln}, we know $\Big( X_T(U_u, \cdot), X_{T/2}(U_u, \cdot)\Big)$ are two correlated Gaussian random variables with mean zero. As long as we can show that $\abs{\text{corr}(X_T(U_u,\cdot), X_{T/2}(U_u,\cdot))} < 1$, then the above probability must be strictly positive. Let us denote 
\begin{align*}
    \mathpzc{f}_t(s) := \frac{1}{\sqrt{L-u}} \Big(e^{-(t-s)\lambda_{-}} - e^{-(t-s)\lambda_{+}}\Big).
\end{align*}
Then 
\begin{align*}
    \ee_{\omega\sim \prob}[X_{T/2}(U_u,\omega)^2] &= \int_{0}^{T/2} \bigl(\mathpzc{f}_{T/2}(s)\bigr)^2\ \ud s;\\
    \ee_{\omega\sim \prob}[X_{T}(U_u,\omega)^2] &= \int_{0}^{T} \bigl(\mathpzc{f}_{T}(s)\bigr)^2\ \ud s;\\
    \ee_{\omega\sim \prob}[X_{T/2}(U_u,\omega) X_{T}(U_u, \omega)] &= \int_{0}^{T/2} \mathpzc{f}_{T/2}(s) \mathpzc{f}_{T}(s)\ \ud s.
\end{align*}
By Cauchy-Schwarz inequality, 
\begin{align*}
   \Abs{\ee_{\omega\sim \prob}\bigl[X_{T/2}(U_u,\omega) X_{T}(U_u, \omega)\bigr]} &=  \Abs{\int_{0}^{T/2} \mathpzc{f}_{T/2}(s) \mathpzc{f}_{T}(s)\ \ud s} \\
   &\le \sqrt{\int_{0}^{T/2} \bigl(\mathpzc{f}_{T/2}(s)\bigr)^2\ \ud s\ \int_{0}^{T/2} \bigl(\mathpzc{f}_{T}(s)\bigr)^2\ \ud s}\\
   &< \sqrt{\int_{0}^{T/2} \bigl(\mathpzc{f}_{T/2}(s)\bigr)^2\ \ud s\ \int_{0}^{T} \bigl(\mathpzc{f}_{T}(s)\bigr)^2\ \ud s} \\
   &= \sqrt{\ee_{\omega\sim\prob}\bigl[X_{T/2}(U_u,\omega)^2\bigr]\ \ee_{\omega\sim\prob}\bigl[X_{T}(U_u, \omega)^2\bigr]}.
\end{align*}
Therefore,
$\abs{\text{corr}(X_{T/2}, X_T)} < 1$, and this completes the proof.
\hfill
\end{proof}

 \end{document}